\documentclass[12pt]{amsart}
\font\emailfont=cmtt10

\usepackage{amsmath,amsthm,amsfonts,amscd,flafter,epsf}
\usepackage{graphics}
\usepackage{epsfig}
\usepackage{psfrag}

\newcommand\commentable[1]{#1}

\newtheorem{theorem}{Theorem}[section]
\newtheorem{prop}[theorem]{Proposition}
\newtheorem{cor}[theorem]{Corollary}

\newtheorem{lemma}[theorem]{Lemma}

\newtheorem{defn}[theorem]{Definition}

\def\endproof{\relax\ifmmode\expandafter\endproofmath\else
  \unskip\nobreak\hfil\penalty50\hskip.75em\hbox{}\nobreak\hfil\bull
  {\parfillskip=0pt \finalhyphendemerits=0 \bigbreak}\fi}
\def\endproofmath$${\eqno\bull$$\bigbreak}
\def\bull{\vbox{\hrule\hbox{\vrule\kern3pt\vbox{\kern6pt}\kern3pt\vrule}\hrule}}

\newcommand{\Q}{\mathbb{Q}}
\newcommand{\R}{\mathbb{R}}

\newcommand{\Z}{\mathbb{Z}}

\newcommand{\OneHalf}{\frac{1}{2}}

\newcommand{\cm}{\cdot}

\newcommand{\F}{\mathbb F}

\newcommand\relspinc{\underline{\spinc}}

\newcommand\Filt{\mathcal F}

\newcommand\x{\mathbf x}

\newcommand\p{\mathbf p}
\newcommand\q{\mathbf q}
\newcommand\y{\mathbf y}

\newcommand\ModFlow{\mathcal M}

\newcommand\ModSphere{\ModFlow\left({\mathbb S}\longrightarrow 
\Sym^{g-1}(\Sigma_{1})\times \Sym^2(\Sigma_{2})\right)}
\newcommand\ModSpheres\ModSphere

\newcommand\CFa{\widehat{CF}}

\newcommand\HFa{\widehat{HF}}

\newcommand\Mas{\mu}
\newcommand\UnparModSp{\widehat \ModSp}
\newcommand\UnparModFlow\UnparModSp
\newcommand\Mod\ModSp

\newcommand{\spinc}{\mathfrak s}
\newcommand{\spincX}{\mathfrak r}

\newcommand\ModMaps{\mathcal M}
\newcommand\ModSp\ModMaps

\newcommand\Ta{{\mathbb T}_{\alpha}}
\newcommand\Tg{{\mathbb T}_{\gamma}}
\newcommand\Tb{{\mathbb T}_{\beta}}
\newcommand\Tc{{\mathbb T}_{\gamma}}

\newcommand\alphas{\mbox{\boldmath$\alpha$}}

\newcommand\betas{\mbox{\boldmath$\beta$}}
\newcommand\gammas{\mbox{\boldmath$\gamma$}}

\newcommand\Triple{$(\Sigma,\alphas,\betas,\gammas,w,z)$}
\newcommand\Dom{\mathcal D}

\newcommand\PerDom{\mathcal P}

\newcommand\Yab{$Y_{\tiny \alpha,\beta}$}
\newcommand\Yag{$Y_{\tiny \alpha,\gamma}$}
\newcommand\Ybg{$Y_{\tiny \beta,\gamma}$}

\newcommand\spincrel\relspinc

\newcommand\CFK{CFK}
\newcommand\HFK{HFK}

\newcommand\CFKa{\widehat\CFK}

\newcommand\CFKinf{\CFK^{\infty}}

\newcommand\HFKa{\widehat\HFK}

\newcommand\BasePt{w}
\newcommand\FiltPt{z}

\newcommand\Dual{\mathcal D}
\newcommand\Duality\Dual

\newcommand\ons{Ozsv{\'a}th and Szab{\'o}}
\newcommand\os{Ozsv{\'a}th-Szab{\'o}}

\newcommand\companion{K}
\newcommand\WDm{$D_+(K,t-1)$}
\newcommand\WDnmm{D_+(K,t-1)}

\newcommand\WD{$D_+(K,t)$}
\newcommand\WDnm{D_+(K,t)}

\newcommand\HFD{$\HFKa(D_+(K,t),1)$} 
\newcommand\CFD{$\CFKa(D_+(K,t),1)$} 
\newcommand\HFDm{$\HFKa(D_+(K,t-1),1)$}

\newcommand\HFDnm{\HFKa(D_+(K,t),1)} 
\newcommand\CFDnm{\CFKa(D_+(K,t),1)} 
\newcommand\HFDnmm{\HFKa(D_+(K,t-1),1)} 
 
\newcommand\HFDgnm{\HFKa_*(D_+(K,t),1)}

\newcommand\HFM{$\underset{\{\spinc_i\in Spin^c(S^3_t(K)\}}\bigoplus \HFKa(S^3_t(K),\mu_K,\spinc_i)$}
\newcommand\CFMnm{\underset{\{\spinc_i\in Spin^c(S^3_t(K)\}}\bigoplus \CFKa(S^3_t(K),\mu_K,\spinc_i)}
\newcommand\HFMnm{\underset{\{\spinc_i\in Spin^c(S^3_t(K)\}}\bigoplus \HFKa(S^3_t(K),\mu_K,\spinc_i)}

\newcommand\Surg{$S^3_t(K)$}
\newcommand\Surgnm{S^3_t(K)}

\newcommand\HFMnospinc{$\bigoplus \HFKa(\Surgnm,\mu_\companion,\spinc_i)$}
\newcommand\WDneg{$D_+(K,-t)$}
\newcommand\WDnmneg{D_+(K,-t)}

\newcommand\HFDgnmneg{\HFKa_*(D_+(K,-t),1)}

\newcommand\Long{$\lambda_P\#\lambda_K$}

\newcommand\Merid{$\mu_P\#\mu_K$}

\commentable{

\title[{Knot Floer homology of Whitehead doubles}] 
{Knot Floer homology of Whitehead doubles}

}

\author[Matthew Hedden]{Matthew Hedden}
\address{Department of
Mathematics, Princeton University\newline
\indent{\emailfont{mhedden@math.princeton.edu}}}
\thanks{The author was supported by an NSF postdoctoral fellowship, and partially supported by the NSF holomorphic curves FRG grant during the course of this work.}

\begin{document}

\begin{abstract}
	In this paper we study the knot Floer homology invariants of the twisted and untwisted Whitehead doubles of an arbitrary knot $K$. We present a formula for the filtered chain homotopy type of $\HFKa(D_{\pm}(K,t))$ in terms of the invariants for $K$, where $D_{\pm}(K,t)$ denotes the $t$-twisted positive (resp. negative)-clasped Whitehead double of $K$. In particular, the formula can be used iteratively and can be used to compute the Floer homology of manifolds obtained by surgery on Whitehead doubles. An immediate corollary is that $\tau(D_+(K,t))=1$ if $t< 2\tau(K)$ and zero otherwise, where $\tau$ is the Ozsv{\'a}th-Szab{\'o} concordance invariant.  It follows that the iterated untwisted Whitehead doubles of a knot satisfying $\tau(K)>0$ are not smoothly slice.   

\end{abstract}

\maketitle
\section{Introduction}
\label{sec:intro}
Satellite knots are frequently studied objects in the world of low-dimensional topology. Among the most famous satellite knots are the Whitehead doubles, which have been at the heart of many beautiful constructions \cite{Akbulut1,Akbulut2,Gompf2,Freedman2}.  As discussed below, the untwisted double of an arbitrary knot has classical invariants such as the Alexander polynomial and signature identical to those of the unknot.  Thus computing values for Whitehead doubles provides an interesting test of any new knot invariant's strength.  Perhaps the Whitehead doubles have shone most brightly  in the study of knot concordance, where they provide examples of knots which are topologically slice yet not smoothly slice.  In this way the Whitehead doubles showcase the remarkable distinction between the smooth and topological categories in dimension four.  

In recent years \ons \ have constructed a comprehensive and powerful set of invariants for low-dimensional topological and geometric objects using the Floer homology theory of pseudo-holomorphic curves.  The purpose of this paper is to study the knot invariants introduced in \cite{Knots,Ras1} in the context of Whitehead doubling.  Our motivation is twofold: $(1)$ to obtain a better understanding of the $(2+1)$ dimensional  topological quantum field theoretic properties of the \os \ invariants and $(2)$ to exploit the power of the invariants to answer topological questions, particularly questions related to smooth knot concordance. 

\begin{figure}
\begin{center}
\psfrag{N}{t}
\psfrag{N+3}{\small{t+3}}
\psfrag{V}{V}
\psfrag{P}{P}
\psfrag{K}{$\nu K$}
\psfrag{D}{\WD}
\psfrag{1}{1}
\psfrag{f}{f}
\psfrag{+}{+}
\psfrag{=}{=}
\caption{\label{fig:Whiteheaddouble}
The positive $t$-twisted Whitehead double, \WD, of the left-handed trefoil.  Start with a twist knot, $P$ with $t$ full twists embedded in a solid torus, $V$.  The $``+"$ indicates the parity of the clasp of $P$. $f$ identifies $V$ with the neighborhood of $K$, $\nu K$, in such a way that the longitude for $V$ is identified with the Seifert framing of $K$.  The image of $P$ under this identification is \WD. The $3$ extra full twists in the projection of \WD \ shown arise from the writhe of the trefoil, $-3$. }
\includegraphics{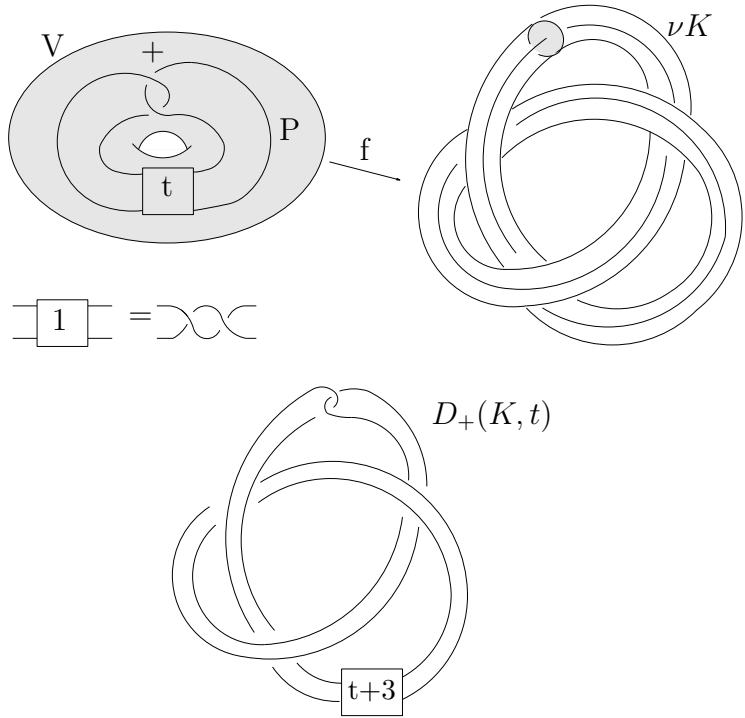}
\end{center}
\end{figure}    
Suppose we have a knot $P$ embedded in a solid torus, $V$.  Letting $K$ be an arbitrary knot, we can identify a tubular neighborhood of $K$ with $V$ in such a way that the longitude of $V$ (the curve on $\partial V$ generating $H_1(V,\Z)\cong\Z$) is identified with a longitude of $K$ coming from a Seifert surface.  The image of $P$ under this identification is a knot, $S$, called a {\em satellite} of $K$.  The knot $P$ is called the {\em pattern} for $S$, while $K$ is referred to as the {\em companion}.  In this language, the {\bf positive $t$-twisted Whitehead double} of a knot $K$ - denoted $D_+(K,t)$ - is a satellite of $K$, where the pattern is a positive-clasped twist knot with $t$ twists.  See Figure \ref{fig:Whiteheaddouble} for an illustration.  Whitehead doubling in the context of \os \ homology was first studied by Eftekhary in \cite{Eaman}.  For other results regarding knot Floer homology and satellite knots, see \cite{Cable,MyThesis,STau,Ording,Ni}.  Before stating the main theorem, recall that associated to an integer homology three-sphere, $Y^3$, is the \os \ chain complex, denoted $\CFa(Y^3)$ (see \cite{HolDisk} for definitions and generalizations).  \ons \ showed that the homology of this chain complex is an invariant of the smooth three-manifold.  In \cite{Knots,Ras1}, it was shown that a knot $K\subset Y^3$ induces a filtration of $\CFa(Y^3)$, and that the filtered chain homotopy type of the resulting filtered chain complex is an invariant of the knot $(Y,K)$.  In the case where $Y^3=S^3$, the three-dimensional sphere, $K$ is a knot in the classical sense and the filtration of $\CFa(S^3)$ is denoted $\Filt(K)$. More explicitly, we have the following increasing sequence of subcomplexes: 
$$ 0=\Filt(K,-i) \subseteq \Filt(K,-i+1)\subseteq \ldots \subseteq \Filt(K,n)=\CFa(S^3),$$
\noindent  We denote the quotient complexes $\frac{\Filt(K,j)}{\Filt(K,j-1)}:=\CFKa(K,j)$, and the homology of these quotients, denoted $\HFKa(K,j)$, are commonly referred to as the {\em knot Floer homology groups} of $K$.  It follows from the fact that the filtered chain homotopy type of $\Filt(K)$ is an invariant of $K$ that the knot Floer homology groups are also knot invariants.  The following theorem suggests that the knot Floer homology groups can be viewed as a ``categorification'' of the symmetrized Alexander-Conway polynomial, in the same spirit that the Khovanov homology groups \cite{Khovanov} are a categorification of the Jones polynomial:

\begin{theorem}(\os \ \cite{Knots}, Rasmussen \cite{Ras1}) Let $K\subset S^3$ be a knot and $\Delta_K(T)$ its Alexander-Conway polynomial.  Then 
$$\sum_{i} \chi \left(\HFKa(K,i)\right) \cm T^i =\Delta_K(T).$$
\end{theorem}

Suppose we look at a satellite $S$ of a knot, $K$, where the pattern $P$ in the construction represents $p$ times a generator of $H_1(V,\Z)$. We have the following classical formula for the Alexander polynomial of $S$ (see \cite{Lickorish}),
$$
\Delta_{S}(T)=\Delta_{P}(T)\cm\Delta_K(T^p).
$$

Since the twist knots used as pattern for the Whitehead double construction represent zero in the first homology of the solid torus,  we see that the Alexander polynomial forgets the knot which we are doubling.  Indeed, the Alexander polynomial of \WD \ is given by $$\Delta_{\WDnm}(T)=-t\cm T + (2t+1) - t\cm T^{-1},$$
\noindent independent of $K$.  In particular, the Alexander polynomial of the $0$-twisted Whitehead double of $K$ is trivial.  It is thus an interesting question to ask how, if at all, the knot Floer homology of \WD \ remembers the knot $K$.  In order to state our theorem, we remark that the knot Floer homology groups $\bigoplus\HFKa(K,i)$ themselves have the structure of a filtered chain complex, endowed with a differential induced from the differential on $\CFKa$.  For a knot of Seifert genus one, this induced  differential decomposes as a sum of three homomorphisms: 

$$d^1_1: \HFKa_*(K,1)\longrightarrow \HFKa_{*-1}(K,0),$$
$$d^0_1: \HFKa_*(K,0)\longrightarrow \HFKa_{*-1}(K,-1),$$
$$d_2: \HFKa_*(K,1)\longrightarrow \HFKa_{*-1}(K,-1).$$
\noindent (See Section \ref{sec:middle} for more details)

\begin{theorem}
	\label{thm:main}
	Let $K\subset S^3$ be a knot with Seifert genus $g(K)=g$.  Then for $t\ge 2\tau(K)$ we have:

	$$\HFKa_*(\WDnm,i)\cong \left\{\begin{array}{ll} 
		\Z_{(1)}^{t-2g-2}  \ \ \bigoplus_{i=-g}^{g}[H_{*-1}(\Filt(K,i))]^2 & i=1\\
		\Z_{(0)}^{2t-4g-3} \ \bigoplus_{i=-g}^{g}[H_{*}(\Filt(K,i))]^4 & i=0\\
		\Z_{(-1)}^{t-2g-2} \ \ \bigoplus_{i=-g}^{g}[H_{*+1}(\Filt(K,i))]^2 & i=-1\\
	\end{array}\right.$$	

Whereas for $t< 2\tau(K)$ the following holds:

	$$\HFKa_*(\WDnm,i)\cong \left\{\begin{array}{ll} 
		\Z_{(1)}^{2\tau(K)-2g-2} \oplus\Z_{(0)}^{2\tau(K)-t} \bigoplus_{i=-g}^{g}[H_{*-1}(\Filt(K,i))]^2 & i=1\\
		\Z_{(0)}^{4\tau(K)-4g-4} \oplus\Z_{(-1)}^{4\tau(K)-2t-1}\bigoplus_{i=-g}^{g}[H_{*}(\Filt(K,i))]^4 & i=0\\
		\Z_{(-1)}^{2\tau(K)-2g-2}  \oplus\Z_{(-2)}^{2\tau(K)-t}  \bigoplus_{i=-g}^{g}[H_{*+1}(\Filt(K,i))]^2 & i=-1\\
	\end{array}\right.$$
Furthermore, $d_2=0$, regardless of $t$, and this together with the formulas above determine the filtered chain homotopy type of $\Filt(\WDnm)$. 	
\end{theorem}

\noindent {\bf Remarks:} The precise way that $d_2=0$ and our formula determine $\Filt(\WDnm)$ is discussed in Section \ref{sec:middle}.  The term $\tau(K)$ is the \os \ concordance invariant \cite{FourBall,Ras1} discussed below.   Since $\tau(K)\le g(K)$, the reader may therefore question what is meant by a term such as $\Z_{(1)}^{2\tau(K)-2g-2}$with negative exponent.  By $\Z^{-n}_{(1)}$, for instance, we mean the quotient of the remaining group by a subgroup of dimension $n$, supported in homological grading $1$.  

Letting $\overline{K}$ denote the reflection of a knot $K$ (i.e. in a given projection for $K$,  $\overline{K}$ is obtained from $K$ by changing each over-crossing to an under-crossing), we have the following formula for the knot Floer homology (Proposition $3.7$ of \cite{Knots}),
\begin{equation}
	\label{eq:mirror}
	\HFKa_*(\overline{K},i)\cong \HFKa_{-*}(K,-i).
\end{equation}	
In light of the following equality,
$$\overline{D_+(\overline{K},-t)}=D_-(K,t)$$
\noindent we see that Theorem \ref{thm:main} yields the complete answer for the negative clasped doubles as well.

One should compare Theorem \ref{thm:main} with the results of \cite{STau} and \cite{Eaman}. Proposition $2.6$ of \cite{STau} computes the Floer homology of a specific Whitehead double of the $(2,n)$ torus knot 
while \cite{Eaman} equates a particular knot Floer homology group of the $0$-twisted Whitehead double with another invariant, the {\em longitude Floer homology}.  Theorem \ref{thm:main} is a significant improvement over either of these results and over any other results concerning the Floer homology of satellite knots.  In fact, the above theorem is a complete answer to the question of Whitehead doubling:  it handles all values of the twisting parameter, $t$, and all the Floer homology groups.  Moreover, we are able to use our formula iteratively. It is interesting to note that Theorem \ref{thm:main} indicates that the direct sum of all the longitude Floer homology groups of $K$ are determined by the Floer homology of $K$, something not clear from their definition.

\subsection{Concordance invariants}
Whitehead doubles have played an interesting role in the study of knot concordance, where they highlight the distinction between the smooth and topological four-ball genus. Moreover, several fundamental remaining open questions in the field of four-dimensional topological surgery are equivalent to questions related to Whitehead doubling \cite{Freedman2,Freedman3}. 

Before going further, recall that the smooth {\em four-ball genus}, $g_4(K)$  of a knot $K$ is the minimum genus of any smoothly properly embedded surface, $(F,\partial F)$, in the four-dimensional ball whose restriction to $\partial F$ is the given knot in $S^3$.
A knot is said to be smoothly {\em slice} if its smooth four-ball genus is zero.  Two knots, $K_1,K_2$ are said to be smoothly {\em concordant} if $K_1\# -\overline{K}_2$ is smoothly slice. Here $-K$ denotes the knot $K$ with reversed orientation. It can be shown that concordance is an equivalence relation on the set of knots and that under this equivalence the set of knots has the structure of an abelian group, the group operation being connected sum, $K_1\# K_2$. We denote this group, the {\em smooth concordance group of knots}, by $\mathcal{C}$.  We can repeat all the above definitions in the topological category, where we require surfaces to be topologically locally flatly embedded. In this case we denote the (topological) four-ball genus and concordance group by $g_4^{top}(K)$ and $\mathcal{C}_{top}$, respectively.

Whitehead doubling is an easy way to produce non-trivial topologically slice knots, as indicated by the following theorem:
\begin{theorem}(Freedman \cite{Freedman2}) Let $K\subset S^3$ be knot which satisfies $\Delta_K(T)=1$.  Then $K$ is topologically slice. That is, $g_4^{top}(K)=0$.
\end{theorem}
\noindent As mentioned above, the $0$-twisted Whitehead double of any knot satisfies $\Delta_{D_\pm(K,0)}(T)=1$, and hence these knots are topologically slice.  While it is easy to see that the Whitehead double of a smoothly slice knot is also smoothly slice, it was shown by several authors that many Whitehead doubles are not smoothly slice \cite{Akbulut1,Gompf2,Rudolph2}.  The existence of a topologically slice knot which is not smoothly slice is interesting in its own right, as it implies the existence of an exotic smooth structure on $\R^4$ (see \cite{Gompf} for a proof).  It is an open question (see Problem $1.38$ on Kirby's list \cite{Kirby}) whether the Whitehead double of $K$ is smoothly slice only when $K$ is smoothly slice.

From the knot Floer homology filtration, we can produce an integer-valued knot invariant $\tau(K)$ useful for the study of smooth knot concordance.  To define it, 
recall that the Floer homology of the three-sphere is isomorphic to $\Z$, supported in homological grading zero.  Thus, one can define the following:
$$\tau(K)=\mathrm{min}\{j\in\Z|i_*:H_*(\Filt(K,j))\longrightarrow H_*(\CFa(S^3))\ ~{\text{~is non-trivial}}\}.$$

\noindent  \ons \ \cite{FourBall} and Rasmussen \cite{Ras1} showed that $\tau(K)$ is an invariant of the smooth concordance class of $K$, and that it provides a bound for the smooth $4$-genus of $K$: $$|\tau(K)|\le g_4(K).$$

\noindent Moreover, $\tau(K)$ is additive under connected sum of knots, and hence provides a homomorphism $\mathcal{C}\rightarrow \Z$.  It is important to note that the knot Floer homology groups of $K$ are in general not sufficient to determine $\tau(K)$ since its definition relies on a more detailed knowledge of the knot filtration $\Filt(K)$. 

Theorem \ref{thm:main} indicates that the dependence of the Floer homology of \WD \ on the twisting parameter is determined by $\tau(K)$.  In fact, a key ingredient used to determine the filtered chain homotopy type of $\Filt(\WDnm)$ is the following:  

\begin{theorem}
	\label{thm:tau}
	$$ \tau(\WDnm)= 
\left\{\begin{array}{ll}

0  & {\text{for
 $t\ge 2\tau(K)$}} \\

 1 & {\text{for
 $t<2\tau(K)$}} \\

\end{array}
\right. $$
\end{theorem}

As a corollary, we can determine the iterated $0$-twisted Whitehead doubles which $\tau$ can be used to show are not smoothly slice.  We let $D_+^i(K)$ denote the $i$-th iterated $0$-twisted Whitehead double of $K$ i.e.  $D_+^1(K)=D_+(K,0)$ and $D_+^i(K)=D_+(D_+^{i-1}(K),0)$:
\begin{cor}
  $\tau(D_+^i(K))\ne 0$ if, and only if, $\tau(K)>0$. Hence, if $\tau(K)>0$ then $D_+^i(K)$ is not smoothly slice for every $i$.  
\end{cor}

The above theorem and corollary should be compared with results of Livingston and Naik \cite{Livingston2} which determine $\tau(\WDnm)$ for all $t$ outside a finite interval.  Using the Floer homology of the branched double cover of $D_+(K,0)$, Manolescu and Owens are able to show that $D_+(K,0)$ is not slice if $\tau(K)>0$.  However, they were unable to determine whether iterated doubles were slice since $\tau$ of these knots was unknown except in the cases computed by Livingston and Naik.

We should also remark that in the case where the companion knot is the $(2,n)$ torus knot, the above result follows from \cite{STau}.  Indeed the main purpose of \cite{STau} was to show that $\tau(K)$ does not equal half the Rasmussen concordance invariant, $s(K)$ \cite{Ras2}.  We believe that Whitehead doubles of knots with $\tau(K)\ne0$ will provide further examples of knots with $2\tau(K)\ne s(K)$.  

In a related direction, the results of \cite{Livingston2} and \cite{STau} indicate that there are two invariants associated to a knot:
$$t_s(K)=\mathrm{min}\{t\in\Z|s(D_+(K,t)=0\}.$$
$$t_\tau(K)=\mathrm{min}\{t\in\Z|\tau(D_+(K,t)=0\}.$$
It follows from the fact that $s$ and $\tau$ are smooth concordance invariants that $t_s,t_\tau$ are also invariants of the smooth concordance class of $K$.  However, Theorem \ref{thm:tau} shows that $t_\tau(K)=2\tau(K)$, and hence provides no new information.  Preliminary calculations indicate that this is not the case with $t_s$ and we consider the question of the behavior of $t_s$ an  interesting question.

\subsection{Applications and Examples}
In the final section of the paper we use our formula for a few simple applications.  For both its own interest, and to illuminate our theorem, we first present a closed formula for the Floer homology of the iterated $0$-twisted doubles of the figure eight knot.  We then use our formula in conjunction with a theorem of \ons \ to determine the Floer homology, $\HFa(S^3_{+1}(\WDnm)$, where $S^3_{+1}(\WDnm)$ is the three-manifold obtained by $+1$ Dehn surgery on \WD.  

\bigskip

\noindent{\bf Organization:} The next section of the paper is devoted to finding an efficient Heegaard diagram for Whitehead doubles.  In Section \ref{sec:meridian} we analyze this diagram and prove that a particular Floer homology group of the Whitehead double is isomorphic to the Floer homology of the meridian of $K$, viewed as a knot in the manifold obtained by $t$-surgery on $K$.  Section \ref{sec:largen} computes these groups for sufficiently large values of the twisting parameter, determining \HFD \ for large $|t|$, in the case of the group corresponding to the top filtration level.  We then use the skein exact sequence for knot Floer homology to calculate \HFD \ for the remaining $t$.  In the course of applying the skein sequence, we will determine $\tau(\WDnm)$.  Section \ref{sec:middle} studies the remaining Floer homology group, and the ``higher differentials'' involved in determining the filtered chain homotopy type of $\Filt(\WDnm)$, thus proving Theorem \ref{thm:main}.  The final section of the paper is dedicated to examples and applications of the main theorem.

\bigskip
\noindent{\bf Acknowledgement:} I wish to thank Eaman Eftekhary, Philip Ording, Peter Ozsv{\'a}th, Jacob Rasmussen, and Zoltan Szab{\'o} for interesting conversations.

\section{A Heegaard Diagram for Whitehead doubles}
\label{sec:heegs}
In this section we recall the definition of a compatible Heegaard and introduce an efficient Heegaard diagram for the Whitehead doubles.  We do not review the basics of knot Floer homology (in particular we assume familiarity with the boundary operator, the definition of the knot filtration, etc.)  For an introduction to Heegaard diagrams for knots and computing knot Floer homology from Heegaard diagrams, see Chapter $2$ of \cite{MyThesis}.  

\begin{defn}
\label{def:hd}
A {\em compatible doubly-pointed Heegaard diagram} for a knot $(Y^3,K)$ (or simply a Heegaard diagram for $(Y^3,K)$) is a collection of data 
$$(\Sigma,\{\alpha_1,\ldots,\alpha_g\},\{\beta_1,\ldots,\beta_{g}\},w,z),$$
where  
\begin{itemize}

	\item $\Sigma$ is an oriented surface of genus g, the {\em Heegaard surface},

	\item $\{\alpha_1,\ldots,\alpha_g\}$ are pairwise disjoint, linearly
  independent embedded circles (the {\em $\alpha$ attaching circles}) which specify a handlebody, $U_\alpha$,
  bounded by $\Sigma$, 

	\item $\{\beta_1,\ldots,\beta_{g}\}$ are pairwise disjoint, linearly
  independent embedded circles which specify a handlebody, $U_\beta$,
  bounded by $\Sigma$ such that $U_\alpha\cup_{\Sigma}U_\beta$ is
  diffeomorphic to $Y^3$,

	\item  $K$ is isotopic to the union of two arcs joined along their common endpoints $w$ and $z$.  These arcs, $t_{\alpha}$ and $t_{\beta}$, are 
	properly embedded and parallel to $\Sigma$ in the $\alpha$ and $\beta$-handlebodies, respectively.  

\end{itemize}

\end{defn}

\noindent{\bf Remarks:} This definition is slightly different than what was originally given in \cite{Knots}.  For a discussion of the two definitions and their equivalence, see \cite{MyThesis}.  Note, too, that we are thinking of knots which may not be embedded in the three-sphere, $S^3$.  If we refer to a knot in $S^3$ we will drop $Y^3$ from the notation.  

\subsection{A diagram for Whitehead doubles}
Knot Floer homology is defined in terms of the doubly-pointed Heegaard diagram described above.  Thus, in order to study Whitehead doubles, we first find a compatible diagram.  

We begin by outlining our strategy for producing a diagram for Whitehead doubles (or more general satellites). This construction will be similar that used by Eftekhary in \cite{Eaman}.  We begin with two Heegaard surfaces, one corresponding to the pattern knot and one to the companion.  On each surface we have an extra $\alpha$ attaching curve (that is, for a surface of genus $g$, we have $(g+1)$ $\alpha$ curves).  If we remove the final two $\alpha$ attaching curves, each diagram specifies a manifold with torus boundary - in the case of the pattern knot the manifold is a solid torus, while for the companion it is the knot complement, $S^3-K$.  Furthermore, the final two $\alpha$ curves intersect each other in a unique point and we can think of these curves as a framing (parametrization) of the torus boundary. Forming the connected sum of the two Heegaard diagrams in a neighborhood of the intersection point of the final two $\alpha$ curves will correspond to identifying the boundary tori of the two three-manifolds.  The homeomorphism identifying the tori will depend both on the $\alpha$ curves used in the framing, and how we identify these curves when we form the connected sum of the diagrams.    For the calculations of knot Floer homology found later in the text, we will benefit from treating the Heegaard diagrams discussed in this section as actually specifying two different manifolds each.  This is due to the presence of the extra $\alpha$ attaching curve - by deleting one or the other of the final two $\alpha$ curves parametrizing the boundary torus, we obtain a Heegaard diagram for a closed three manifold.  The chain complex for the Floer homology of the Whitehead double discussed in subsequent sections will decompose nicely along the chain complexes associated to the various Heegaard diagrams obtained by using the different $\alpha$ attaching curves. With the idea in place, we begin.

\begin{figure}
\begin{center}
\psfrag{m}{$\beta_2=\mu$}
\psfrag{l}{$\lambda_P$}
\psfrag{mk}{$\mu_P$}
\psfrag{b1}{$\beta_1$}
\psfrag{w}{$z$}
\psfrag{z}{$w$}
\psfrag{y}{$y$}
\psfrag{x}{$x$}
\psfrag{a1}{$a_1$}
\psfrag{a2}{$a_2$}
\psfrag{a3}{$a_3$}
\psfrag{a4}{$a_4$}
\psfrag{alpha1}{$\alpha_1$}
\caption{\label{fig:Pattern}
Genus $2$ Heegaard diagram for the pattern knot in the Whitehead double construction. It is actually two Heegaard diagrams, as described in the text, depending on whether we use the $\alpha$ curve $\lambda_P$ or $\mu_P$.  Note the large black disc where $\lambda_P$ intersects $\mu_P$ - we will glue the Heegaard diagram for the companion knot to this Heegaard diagram along the black disc. }
\includegraphics{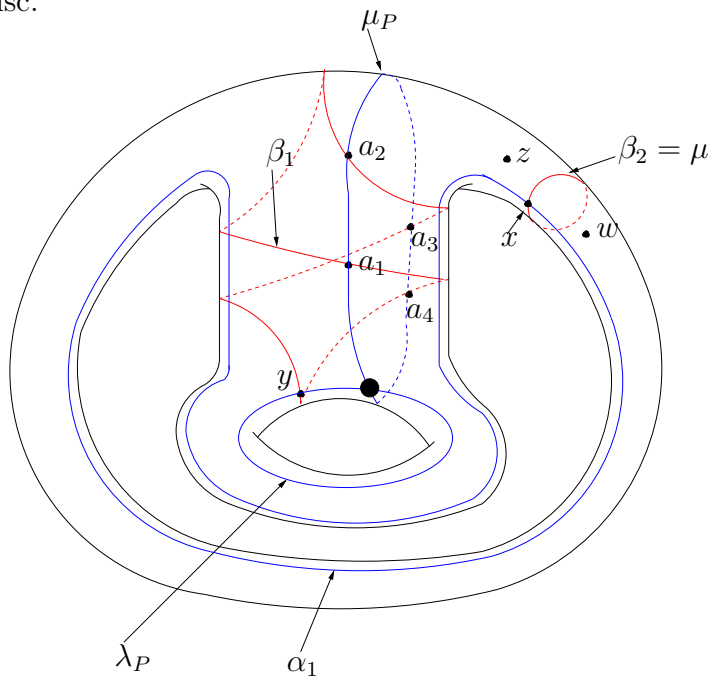}
\end{center}
\end{figure}

Figure \ref{fig:Pattern} depicts a Heegaard diagram associated to the pattern knot, $P$. As mentioned above, there are too many $\alpha$ attaching curves to specify a three-manifold.  We interpret the diagram as two diagrams:
$$\mathrm{hd}(P)=(\Sigma_2,\{\alpha_1,\alpha_2=\lambda_P\},\{\beta_1,\beta_2=\mu\},w,z)$$
$$\mathrm{hd}(Hopf)=(\Sigma_2,\{\alpha_1,\alpha_2=\mu_P\},\{\beta_1,\mu\},w,z).$$

\begin{figure}
\begin{center}
	\psfrag{0}{$0$}
\psfrag{H}{$Hopf$}
\psfrag{mk}{$\mu_K$}
\psfrag{t}{$t$}
\caption{\label{fig:Hopf} The ``knotification'' of the two component positive Hopf link, and the knot $(\Surgnm,\mu_K)$}
\includegraphics{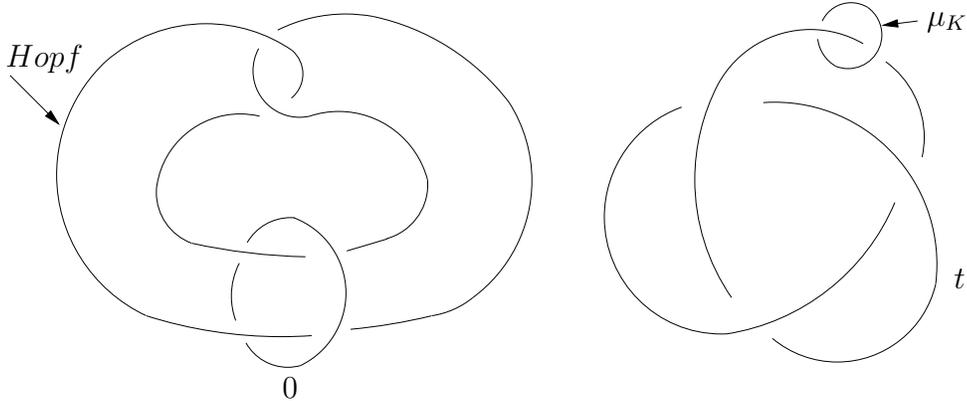}
\end{center}
\end{figure}

Note that $\mathrm{hd}(P)$ specifies the pattern knot in $S^3$, while the diagram $\mathrm{hd}(Hopf)$ specifies the knot in $S^1\times S^2$ shown in Figure \ref{fig:Hopf}. Note also that $\lambda_P \cap \mu_P = \{1 \ \mathrm{ point}\}$, which we draw on the diagram as a black hole.  (The terminology $\mathrm{hd}(Hopf)$ is explained as follows: In \cite{Knots} \ons \ describe a way to associate a null-homologous knot $(\#^{|L|-1}S^1 \times S^2,k(L))$ to a link $(S^3,L)$ of $|L|$ components. Our notation is explained by the fact that the Heegaard diagram specifies the knot $(S^1\times S^2,k(L))$ associated to the two-component Hopf link.)  

\begin{figure}
\begin{center}
\psfrag{m}{$\mu_K$}
\psfrag{l}{$\lambda_K$}
\psfrag{p}{$p$}
\psfrag{b1}{$\beta_j'$}
\psfrag{w}{$w'$}
\psfrag{z}{$z'$}
\psfrag{z''}{$z''$}
\psfrag{y}{$y$}
\psfrag{x}{$x$}
\psfrag{ldots}{$\ldots$}
\psfrag{a}{$\alpha_i'$}
\psfrag{b2}{$\beta_k'$}
\psfrag{Sigma}{$\Sigma_g$}
\psfrag{g}{$p$}
\caption{\label{fig:Companion}
Heegaard diagram for an arbitrary companion knot, $\companion$, in the Whitehead double construction.  We show only the last segment of the diagram which includes a meridian for $\companion$.  As before, it is actually two Heegaard diagrams depending on whether we use the $\alpha$ curve $\lambda_\companion$ or $\mu_\companion$.  Again we denote with a large black disc the intersection of $\lambda_\companion$ and $\mu_\companion$.  }
\includegraphics{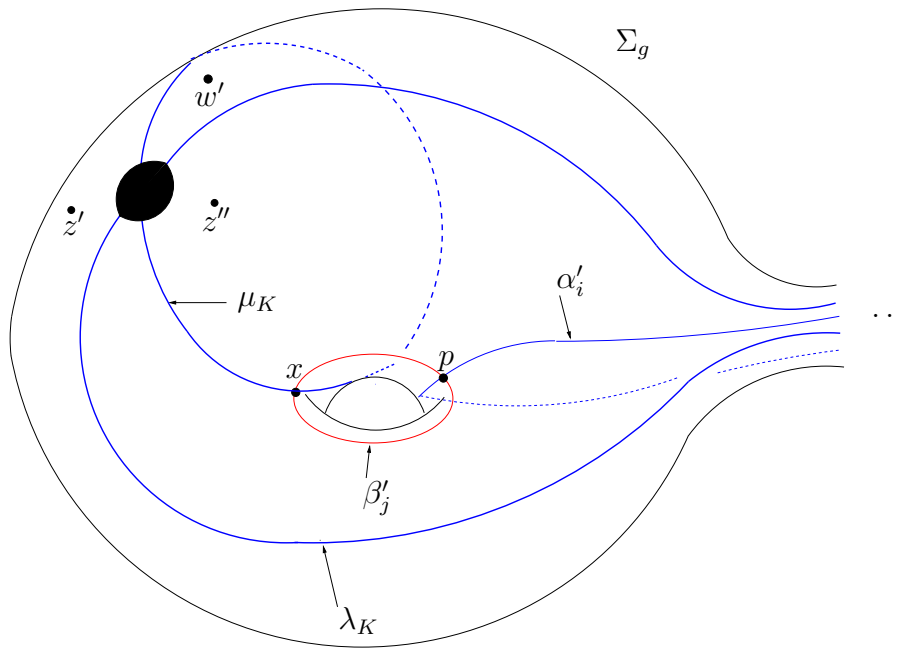}
\end{center}
\end{figure}

We now consider a Heegaard diagram for the companion knot, $\companion$.  In addition to the requirements of Definition \ref{def:hd}, for this diagram we require that one of the $\alpha$ attaching curves, $\mu_\companion$, is a meridian for $\companion$ (so that the diagram without $\mu_\companion$ specifies the knot complement $S^3-\companion$).  This added requirement allows us to draw a framed longitude for $\companion$ embedded on the Heegaard surface as follows: connect $z'$ to $w'$ by a small arc,  $t_{\alpha}$, which intersects only $\mu_\companion$ and an arc $t_{\beta}$ which only intersects the $\alpha$ curves.  The union $\lambda_\companion=t_\alpha \cup_{\{z',w'\}} t_\beta$ is a longitude for the companion knot.  It is clearly embedded on the surface and hence we can think of it as an attaching curve for the $\alpha$-handlebody.   With this extra curve, the genus $g$ diagram has too many $\alpha$ curves, and so, as above, we view it as two diagrams, see Figure \ref{fig:Companion}):
$$\mathrm{hd}(\companion)=(\Sigma_g,\{\alpha_1',\ldots,\alpha_g'=\mu_\companion\},\{\beta_1',\ldots,\beta_{g}'\},w',z').$$
$$\mathrm{hd}( (S^3_t(\companion),\mu_\companion) )=(\Sigma_g,\{\alpha_1',\ldots,\alpha_g'=\lambda_\companion\},\{\beta_1',\ldots,\beta_{g}'\},w',z'').$$

\noindent $\mathrm{hd}(\companion)$ is simply a diagram for the companion knot in $S^3$ with the last $\alpha$ curve a meridian for $\companion$.  The second diagram no longer specifies $S^3$, but instead the manifold obtained by $t$-surgery on the companion knot, $S^3_t(\companion)$, where $t$ is the framing of the longitude $\lambda_\companion$.  We can vary the framing by letting $\lambda_\companion$ circle more or fewer times around the meridian, but $\lambda_\companion \cap \mu_\companion = \{1\ \mathrm{point}\}$, regardless of the framing.

The notation $\mathrm{hd}( (S^3_t(\companion),\mu_\companion) )$ is explained by the fact that the knot in $S^3_t(\companion)$ determined by $w'$ and $z''$ is the meridian of $\companion$. To see this, simply connect $w'$ to $z''$ by arcs in the $\alpha$ and $\beta$ handlebodies for $\mathrm{hd}( (S^3_t(\companion),\mu_\companion) )$ to recover a knot isotopic to $\mu_\companion$.  It should be noted that $(S^3_t(\companion),\mu_\companion)$ is {\em not} a null-homologous knot.  However, when $t\ne0$, $(S^3_t(\companion),\mu)$ is rationally null-homologous i.e. $[\mu]=0\in H_1(S^3_t(\companion),\Q)$.  In this case \ons \ have defined knot Floer homology groups associated to $(S^3_t(\companion),\mu_\companion)$, see \cite{RationalSurgeries}.  
In the case $t=0$ the knot Floer homology groups associated to $(S^3_t(\companion),\mu_\companion)$ are what Eftekhary refers to as the longitude Floer homology groups of $\companion$, see \cite{Eaman}.

Out of the four Heegaard diagrams $\mathrm{hd}(P),\mathrm{hd}(Hopf),\mathrm{hd}(\companion),\mathrm{hd}( (S^3_t(\companion),\mu) )$ we form a single diagram for the $t$-twisted Whitehead double \WD \ by the following construction, which can be found in \cite{Eaman}.

\begin{figure}
\begin{center}

\psfrag{m}{$\mu_P \# \mu_\companion$}
\psfrag{l}{$\lambda_P \# \lambda_\companion$}
\psfrag{Sigma}{$\Sigma_2\# \Sigma_g$}

\psfrag{m1}{$\mu$}
\psfrag{z}{$w$}

\psfrag{w}{$z$}

\caption{\label{fig:Gluing2} Depiction of the gluing process used to obtain a Heegaard diagram for the Whitehead double of $\companion$. We form the connected sum of the diagrams for $P$ and $\companion$ along the black discs, with the diagram for $P$ (Figure \ref{fig:Pattern}) embedded in the $\alpha$ handlebody for $\companion$ (Figure \ref{fig:Companion}).  We illustrate the $\alpha$ curves $\mu_P\# \mu_\companion$ and $\lambda_P\# \lambda_\companion$ which we obtain by connecting the curves on the two original Heegaard diagrams. The only other attaching curve shown is the meridian of the Whitehead double, $\mu$.  The rest of the $\alpha$ and $\beta$ attaching curves are inherited without modification from the diagrams for $P$ and $\companion$.}

\includegraphics{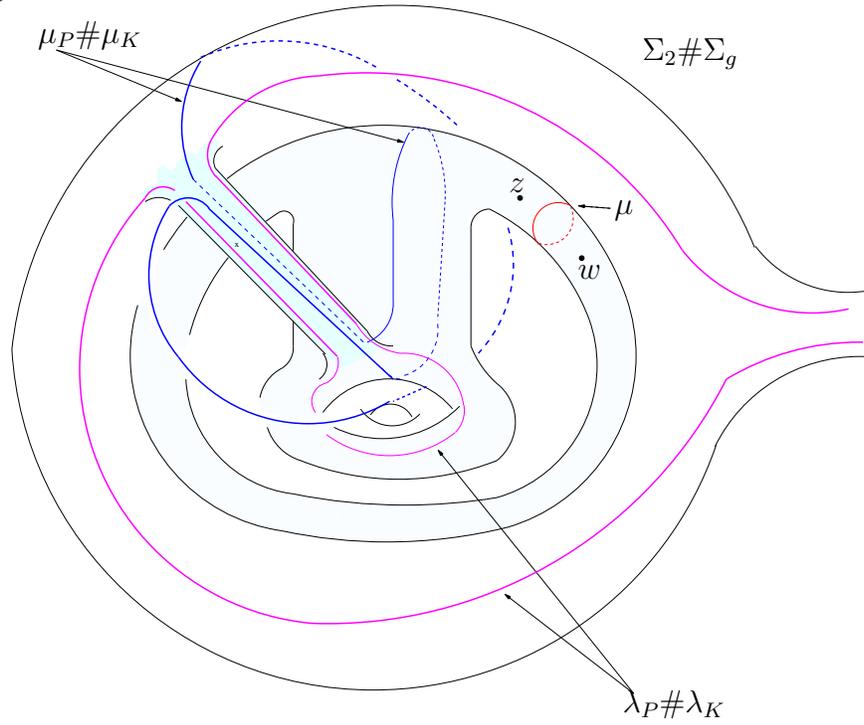}
\end{center}
\end{figure}

We first describe the surface.  Start by embedding the genus two surface for $P$ inside the $\alpha$ handlebody specified by $\mathrm{hd}(\companion)$.  This is shown in Figure \ref{fig:Gluing2}.

Next form the connected sum of the outside surface with the inside surface.  We form the connected sum in a neighborhood of the intersection points  $\lambda_P\cap\mu_P$ and $\lambda_\companion\cap\mu_\companion$, respectively. The resulting surface has genus $g+2$, where $g$ is the genus of the diagram for $\companion$.  Next, we specify the attaching curves.  The $\beta$ attaching curves will be exactly the $\beta$ curves present on the original diagrams. The $\alpha$ curves will also be those present on the original diagrams, except that we connect the longitude and meridian curves, i.e.  $\lambda_P\#\lambda_\companion$ and $\mu_P \# \mu_\companion$. We do this so that the attaching disks for $\lambda_P\#\lambda_\companion$ and $\mu_P \# \mu_\companion$ will ``fill in'' the space between the boundary tori associated to the solid torus $V$ and to $S^3-K$ which is left after forming the connect sum of $\Sigma_g$ and $\Sigma_2$.   Finally, for the basepoints we use the points $z$ and $w$ from the pattern.  Summarizing, we have a diagram: 
$$\mathrm{hd}(\WDnm)=(\Sigma_{g+2},\{\alpha_1,\alpha_1'\ldots,\alpha_{g-1}',\lambda_\companion\#\lambda_\companion,\mu_\companion\#\mu_\companion \},\{\beta_1,\mu,\beta_1'\ldots,\beta_{g}'\},w,z).$$

\begin{figure}
\begin{center}

\psfrag{m}{$\mu$}
\psfrag{l}{$\lambda_P \# \lambda_\companion$}

\psfrag{mk}{$\mu_P \# \mu_\companion$}
\psfrag{A}{A}
\psfrag{B}{B}
\psfrag{C}{C}
\psfrag{D}{D}
\psfrag{E}{E}
\psfrag{F}{F}
\psfrag{b1}{$\beta_1$}
\psfrag{w}{$z$}
\psfrag{z}{$w$}
\psfrag{y}{$y$}
\psfrag{x}{$x$}
\psfrag{a1}{$a_1$}
\psfrag{a2}{$a_2$}
\psfrag{a3}{$a_3$}
\psfrag{a4}{$a_4$}
\psfrag{alpha1}{$\alpha_1$}
\caption{\label{fig:Handleslide} Heegaard moves demonstrating that hd$(\WDnm)$ specifies $S^3$. (A) Begin with hd(\WD).  (B) Handleslide $\beta_1$ over $\mu$.  (C) Perform an isotopy. (D) - Another handleslide of $\beta_1$ over $\mu$ and isotopy of $\beta_1$.  (E) Destabilization (handle-cancellation) of $\mu$ and $\alpha_1$.  (F) Destabilization of $\beta_1$ with $\lambda_P \# \lambda_\companion$.  All Heegaard moves occur in the complement of the basepoint $w$, useful for determining gradings in Section \ref{sec:largen}.}

\includegraphics{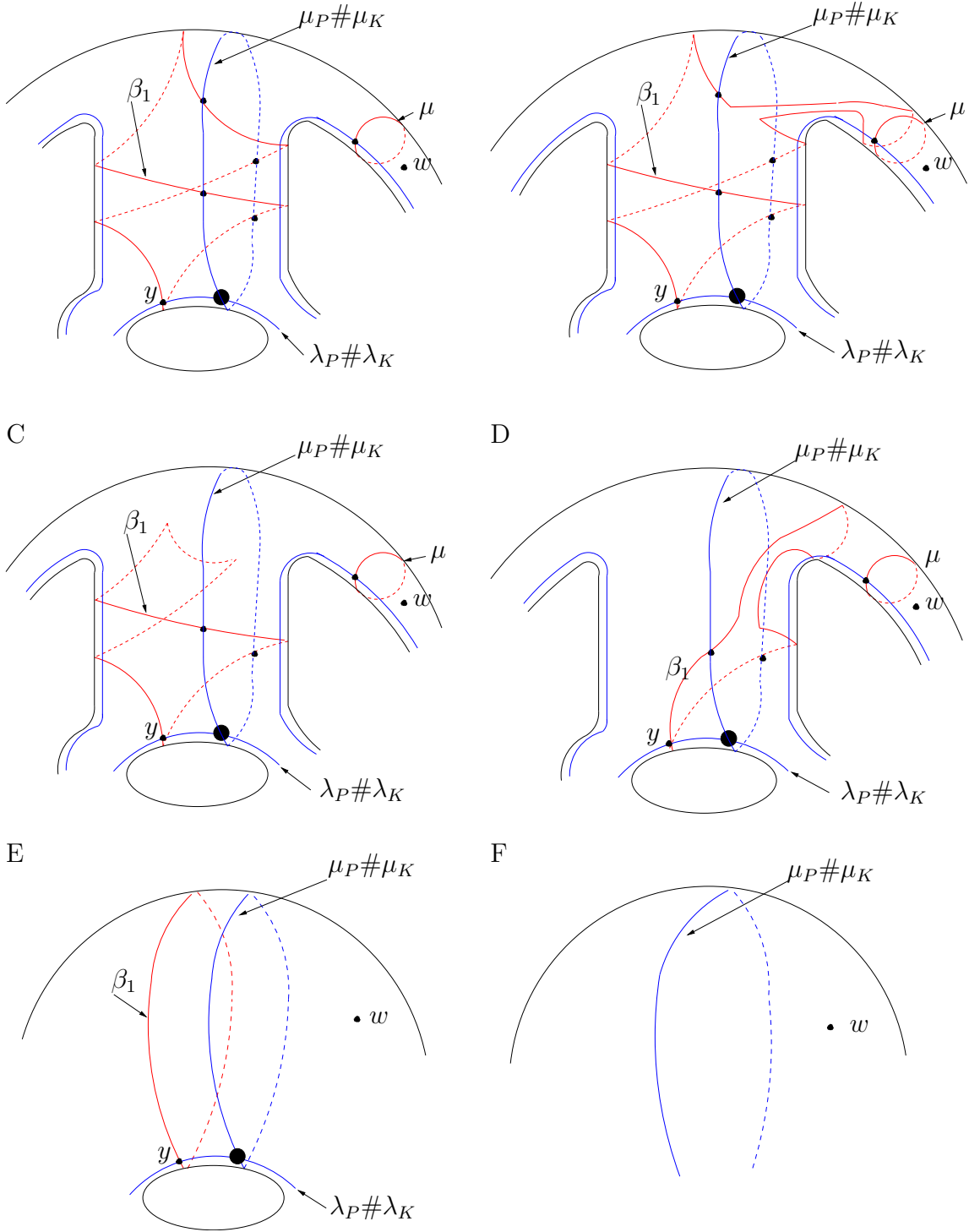}
\end{center}
\end{figure}

We must verify that this diagram is compatible with the $t$-twisted Whitehead double of $P$.  We first show that the three-manifold specified by the diagram is $S^3$.  To see this, first handle-slide $\beta_1$ over $\mu$.  After isotoping $\beta_1$, we handle-slide it over $\mu$ again.  After this second handle-slide, it is immediate that $\mu$ and $\alpha_1$ form a canceling $1-2$ handle pair and can be destabilized, see Figure \ref{fig:Handleslide}.  After destabilizing, $\beta_1$ now forms a canceling pair with the curve $\lambda_\companion \# \lambda_P$.   Destabilizing this pair leaves us with the diagram hd$(\companion)$, which specified $S^3$ by assumption.  

 It remains to see that hd(\WD) specifies the Whitehead double.  However, this can be easily verified by drawing a longitude for the knot specified by hd(\WD) in the same way we drew longitudes for the various diagrams used in the construction.  Indeed, the longitude $\lambda_P$ for the solid torus $V$ in which $P$ was embedded is now isotopic to a $t$-framed longitude for the companion $\lambda_\companion$ via an isotopy along the attaching disk for $\lambda_\companion \# \lambda_P$.   It follows that the knot is isotopic to the $t$-twisted Whitehead double of $\companion$.  This completes our construction of the Heegaard diagram for \WD.  

\noindent {\bf Remark:} Lipshitz \cite{Lipshitz2} is developing a Heegaard Floer invariant for three-manifolds with parametrized boundary.  We can understand the diagrams presented here from his perspective as follows.  By removing a disc from each Heegaard surface in a neighborhood of the intersection point between the final two $\alpha$ curves (the black hole in Figures \ref{fig:Pattern} and \ref{fig:Companion}) we obtain a Heegaard surface with boundary.  The final two $\alpha$ curves parametrizing our boundary then become properly embedded essential arcs on the punctured Heegaard surface.  This is the diagram used by Lipshitz to define his invariant.

\section{Identification of $\widehat{HFK}(D_+(K,t),1)$ with $\widehat{HFK}(S^3_t(K),\mu_K)$}
\label{sec:meridian}
In this section we examine the Heegaard diagram for the Whitehead doubles constructed in the previous section. We will first discuss the generators of the knot Floer homology chain complex and separate them into their respective filtration levels.  Isolating our attention to the top filtration level, we will find a natural identification of the chain complexes: $$ \CFDnm = \CFMnm, $$  
\noindent where the second chain complex is the direct sum of the knot Floer homology chain complexes associated to the meridian of the companion knot, $\mu_\companion$, viewed as a knot in the manifold obtained by $t$-framed Dehn surgery on $K$, $S^3_t(\companion)$.  Note that this is really a double sum taken over Spin$^c$ structures on $S^3_t(\companion)$ and over filtration levels induced by relative Spin$^c$ structures on $\Surgnm-\mu_\companion$.  Let us begin by identifying the generators of the knot Floer homology chain complex associated to the diagram hd(\WD).

 Our first observation is that the generators of the chain complex naturally split into two types, where the splitting is in terms of the four diagrams we used to construct hd(\WD).  The types are of the form:

\begin{enumerate}
	\item $\{x,y\} \times \p \in \CFKa(P) \times \CFKa(\companion)$
	\item $\{x,a_i\} \times \q \in \CFKa(S^1\times S^2, Hopf) \times \CFa(\Surgnm)$
\end{enumerate}

The splitting occurs because the generators correspond to $(g+2)$-tuples of intersection points between the $\alpha$ and $\beta$ attaching curves, with each $\alpha$ and $\beta$ curve used exactly once.  Since \Long \ and \Merid  \ are the only two attaching curves which pass through the connect sum region the splitting is determined by the surface, $\Sigma_2$ or $\Sigma_g$, where the \Long \ component of the $(g+2)$-tuple lies.

\space

\noindent {\bf Note}: Since the intersection point $\{x\}$ occurs as part of the $2$-tuple corresponding to any generator of $\CFKa(P)$ or $\CFKa(S^1\times S^2,Hopf)$, we will subsequently drop it from the notation e.g.  $\{y\}:=\{x,y\}$.

We turn our attention to understanding the relative filtration difference between pairs of generators.  To do this, we identify Whitney disks connecting pairs of generators in the $(g+2)$-fold symmetric product of the Heegaard surface.  In fact, we find it more convenient to identify the {\em domains} of Whitney disks, by which we mean $2$-chains lying in $\Sigma_{g+2}$ with boundary in the attaching curves, and corner points contained in the $(g+2)$-tuple of intersection points representing the generators.  For the equivalence of these methods, see \cite{Ras1,MyThesis}.  Before beginning, we recall the following definition, found in \cite{HolDisk}:

\begin{defn}
	A {\em periodic domain} is a $2$-chain, $\PerDom \subset \Sigma$, such that the boundary of $\PerDom$ consists of a linear combination of attaching circles and so that the local multiplicity at $w$ is zero i.e $\partial \PerDom = \Sigma_{i=1}^g [n_i \alpha_i + m_i \beta_i], n_i,m_i\in \Z$ and $n_w(\PerDom)=0$.

\end{defn}

\begin{figure}
\begin{center}
	\psfrag{Front}{Front}
	\psfrag{Back}{Back}
\psfrag{m}{$\beta_2=\mu$}
\psfrag{l}{$\lambda_P$}
\psfrag{mk}{$\mu_P$}
\psfrag{b1}{$\beta_1$}
\psfrag{w}{$z$}
\psfrag{z}{$w$}
\psfrag{y}{$y$}
\psfrag{x}{$x$}
\psfrag{a1}{$a_1$}
\psfrag{a2}{$a_2$}
\psfrag{a3}{$a_3$}
\psfrag{a4}{$a_4$}
\psfrag{1}{$1$}
\psfrag{-1}{$-1$}
\psfrag{0}{$0$}

\psfrag{alpha1}{$\alpha_1$}
\caption{\label{fig:PerDom}
Illustration of the generator for the space of periodic domains for the diagram hd$(Hopf)$.  Non-zero multiplicities of the two-chain are indicated with shading.  Dark red indicates multiplicity $-1$, while light blue indicates multiplicity $1$. }
\includegraphics{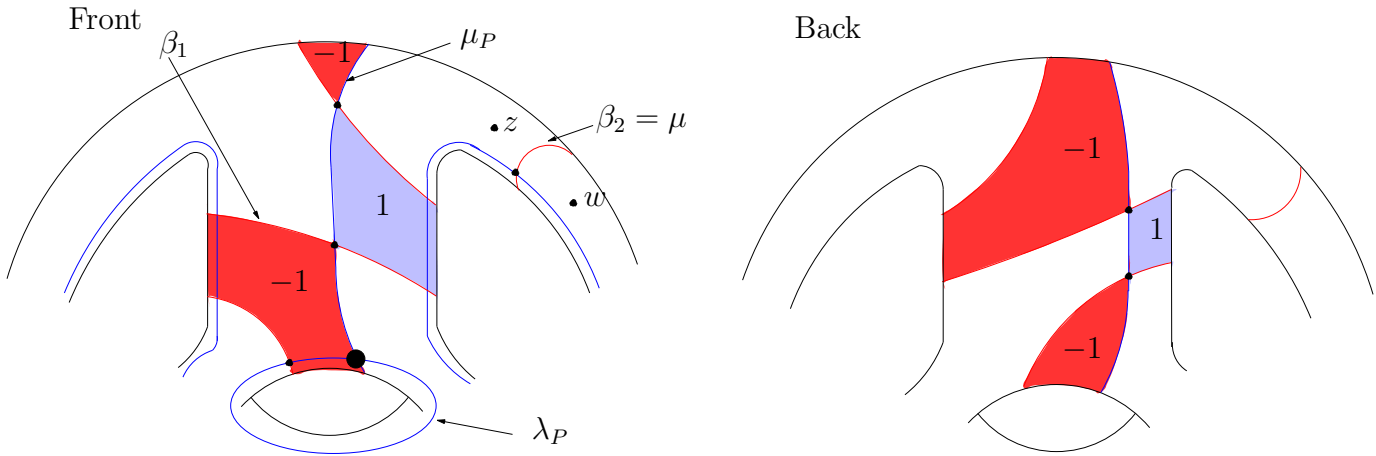}
\end{center}
\end{figure}

If $H_1(Y,\Q)=0$, it follows that there will not exist any periodic domains.  Indeed, since hd(\WD) is a diagram for $S^3$ it will not contain any periodic domains.  However, the diagram hd$(Hopf)$ which went into the construction of hd(\WD) represented $S^1\times S^2$, and this diagram contains period domains.  A generator for the space of periodic domains on hd$(Hopf)$ is shown in Figure \ref{fig:PerDom}.  We will subsequently refer to this generator as $\PerDom$.

To begin, we determine the filtration difference between pairs of points of Type $(1)$.

\begin{lemma}
Suppose $\{y\} \times \p_i,  \{y\} \times \p_k, \in \CFKa(P) \times \CFKa(\companion)$.   Then

$$\Filt(\{y\} \times \p_i)-  \Filt(\{y\}\times \p_k) = 0. $$
 
\end{lemma}

\begin{proof}  We first note that since $\p_i$ and $\p_k$  can be viewed as generators of $\CFKa(\companion)$, they can be connected by a Whitney disk, $\phi$, with domain contained in hd$(\companion)$.  However, if $\p_i$ and $\p_k$ (viewed as generators of $\CFKa(\companion)$) lie in different filtration levels,  then the boundary of  $\phi$ must contain the meridian $\mu_\companion$ with non-zero multiplicity (this follows from the definition of the filtration).  We can complete such a Whitney disk for $\p_i,\p_k\in\CFKa(\companion)$ to a disk connecting $\{y\} \times \p_i$ to $  \{y\}\times \p_k$ in hd(\WD) by forming the boundary sum of $\phi$ with $n\cm \PerDom$ where $n$ is the filtration difference of $\p_i$ and $\p_k$.  The lemma follows because $n_z(\PerDom)=n_w(\PerDom)=0$.  
\end{proof}

Next, we handle the filtration difference of pairs of points of Type $(2)$.

\begin{lemma}

\label{lemma:Lemma2}

Suppose $\{a_i\} \times \q_j,  \{a_i\} \times \q_k, \in \CFKa(S^1\times S^2, Hopf) \times \CFa(\Surgnm)$   Then

$$\Filt(\{a_i\} \times \q_j)-  \Filt(\{a_i\}\times \q_k) = 0. $$
 
\end{lemma}

\begin{proof}  At first sight it appears the method used in the proof of the preceding lemma is hopeless.  Since $H^2(\Surgnm,\Z)\ne0$, we cannot connect arbitrary pairs of  generators for  $\CFa(\Surgnm)$ with a Whitney disk - $\q_j$ and $\q_k$ could represent different Spin$^c$ structures on \Surg.  However, the obstruction to finding a Whitney disk connecting  $\q_j$ to $\q_k$ lies in $H^2(\Surgnm,\Z) \cong H_1(\Surgnm,\Z)$ which is generated by the meridian of $\companion$, $\mu_\companion$.  Thus, if we allow Whitney disks whose boundary contains $\mu_\companion$ (not an attaching curve for hd(\Surg), but present on the diagram for \WD) we can connect $\q_j$ and $\q_k$ regardless of their Spin$^c$ structure.  By completing these Whitney disks with periodic domains as in the proof above, we recover the lemma. 
\end{proof}

Now we examine the effect of varying the generator on the diagram, hd$(Hopf)$.

\begin{lemma} 
	\label{lemma:filt}
	For all $\p \in \CFa(\Surgnm)$ we have: 
$$\Filt(\{a_1\} \times \p) -  \Filt(\{a_4\}\times \p) = 1. $$
$$\Filt(\{a_2\} \times \p) -  \Filt(\{a_1\}\times \p) = -1. $$
$$\Filt(\{a_3\} \times \p) -  \Filt(\{a_4\}\times \p) = -1. $$
$$\Filt(\{a_2\} \times \p) -  \Filt(\{a_3\}\times \p) = 1. $$

Furthermore, there exists a pair $\q \in \CFKa(\companion)$ and $\p \in \CFa(\Surgnm)$ so that: 
$$\Filt(\{y\} \times \q) -  \Filt(\{a_4\}\times \p) = 0. $$ 
\end{lemma}

\begin{figure}
\begin{center}
	\psfrag{Front}{Front}
	\psfrag{Back}{Back}
\psfrag{m}{$\beta_2=\mu$}
\psfrag{l}{$\lambda_P$}
\psfrag{mk}{$\mu_P$}
\psfrag{b1}{$\beta_1$}
\psfrag{w}{$z$}
\psfrag{z}{$w$}
\psfrag{y}{$y$}
\psfrag{x}{$x$}
\psfrag{a1}{$a_1$}
\psfrag{a2}{$a_2$}
\psfrag{a3}{$a_3$}
\psfrag{a4}{$a_4$}
\psfrag{1}{$1$}
\psfrag{-1}{$-1$}
\psfrag{0}{$0$}
\psfrag{alpha1}{$\alpha_1$}
\caption{\label{fig:Disk1}
Illustration of the domain of a Whitney disk connecting $a_2$ to $a_3$ satisfying $n_w(\phi)=0$, $n_z(\phi)=1$. It is topologically an annulus, and can be shown to admit a unique holomorphic representative.  This fact will be used in Section \ref{sec:middle}.  Note the orientation of the Heegaard surface is such that the inward normal forms an oriented basis for $\R^3$.  However, since the ``back'' side of the diagram is actually the mirror image, the orientation of the plane of the page is reversed.  This is our convention on orientations of the pattern knot surface throughout the text.}
\includegraphics{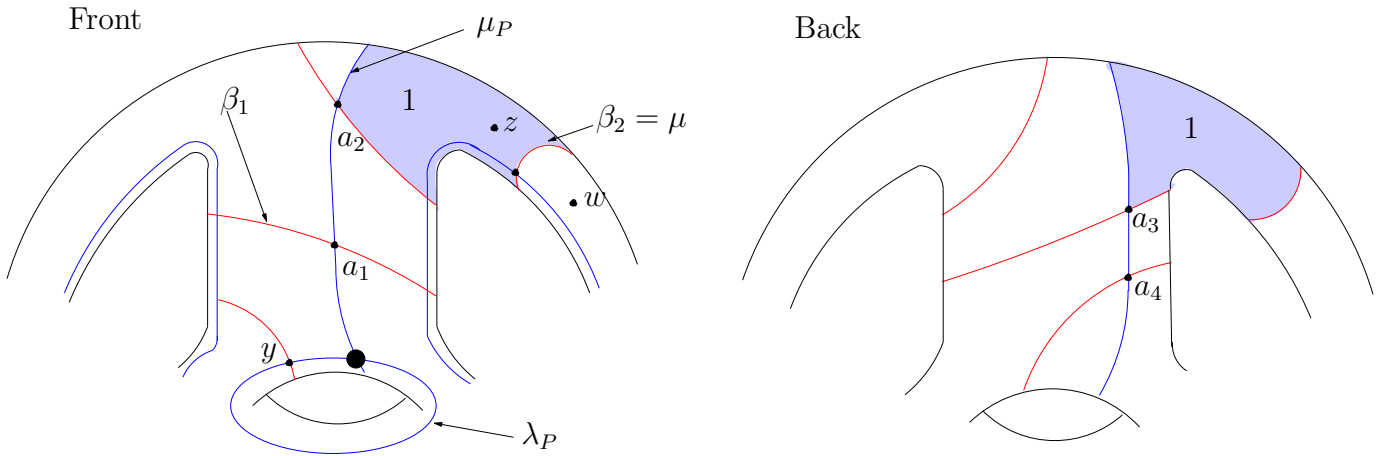}
\end{center}
\end{figure}

\begin{figure}
\begin{center}
	\psfrag{Front}{Front}
	\psfrag{Back}{Back}
\psfrag{m}{$\beta_2=\mu$}
\psfrag{l}{$\lambda_P$}
\psfrag{mk}{$\mu_P$}
\psfrag{b1}{$\beta_1$}
\psfrag{w}{$z$}
\psfrag{z}{$w$}
\psfrag{y}{$y$}
\psfrag{x}{$x$}
\psfrag{a1}{$a_1$}
\psfrag{a2}{$a_2$}
\psfrag{a3}{$a_3$}
\psfrag{a4}{$a_4$}
\psfrag{1}{$1$}
\psfrag{-1}{$-1$}
\psfrag{0}{$0$}
\psfrag{alpha1}{$\alpha_1$}
\caption{\label{fig:Disk2}
Illustration of the domain of a Whitney disk connecting $a_2$ to $a_1$ satisfying $n_w(\phi)=1$, $n_z(\phi)=0$. It is topologically an annulus. }
\includegraphics{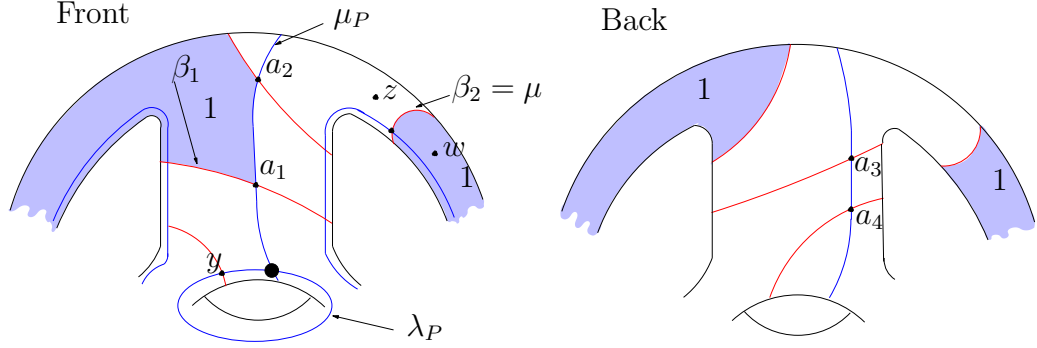}
\end{center}
\end{figure}

\begin{proof} 
	We prove this by explicitly identifying the domains of Whitney disks connecting generators of the above form.  Restricting attention to the Heegaard diagram  hd$(Hopf)$ we can connect pairs of generators $\{a_i\}$ and $\{a_j\}$ with Whitney disks whose domains are topologically annuli.  We complete the annuli to domains on hd(\WD) by extending the annulus on hd$(Hopf)$ by the domain on hd$(\Surgnm)$ which has multiplicity zero everywhere.   See Figures \ref{fig:Disk1} and \ref{fig:Disk2} for illustrations. The first four rows in the table below describe these annuli. The third column lists the Maslov index of the disk (we will need this when we discuss gradings in the next section).  The last two columns indicate the multiplicities of the domains at the points $z$ and $w$.   This completes the first four parts of the lemma.  For the last, we can connect  $\{y\} \times \q $ to  $\{a_4\}\times \p$ by a  domain whose multiplicities are shown in the table. The topology of this domain depends on the framing framing, $t$, and we must choose $\q$ carefully in order to construct it.  This is explained in detail in Subsection \ref{subsec:gradings}.  When $t>0$, Figure \ref{fig:positiverectangle} in Subsection \ref{subsec:gradings} illustrates a domain connecting $\{a_4\}\times \p$ to $\{y\} \times \q $.  When $t<0$ there is a disk $\phi$ connecting$\{a_1\}\times \p$  to $\{y\} \times \q $ with $n_z(\phi)=1$. In either event, we have that
$$\Filt(\{y\} \times \q) -  \Filt(\{a_4\}\times \p) = 0, $$
	\noindent as claimed.

\bigskip
\begin{center}
\begin{tabular}{|c|c||c|c|c|}
\noalign{\hrule height0.8pt}
Start Pt  & End Pt  & $\mu(\phi)$ & $n_{z}(\phi)$ & $n_w(\phi)$  \\
\hline
$\{a_1\}$  & $\{a_4\}$    & $1$ & $1$ & $0$\\
\hline
$\{a_2\}$  & $\{a_1\}$  & $1$ & $0$ & $1$ \\
\hline
$\{a_3\}$  & $\{a_4\}$  & $1$ & $0$ & $1$ \\
\hline
$\{a_2\}$  & $\{a_3\}$  & $1$ & $1$ & $0$  \\
\hline
$\{a_4\}$  & $\{y\}$  & $1$ & $0$ & $0$\\
\hline
\noalign{\hrule height 0.8pt}
\end{tabular}
\end{center}
\end{proof}

The above lemmas are enough to determine the relative filtration difference between any two generators and we see that the chain complex for \WD \  splits into three distinct filtration levels.   The filtration levels have the following form:
$$\CFKa(\WDnm,1)\sim \{a_1\}\times \CFa(\Surgnm)$$
$$\CFKa(\WDnm,0)\sim [\{a_2\}\times \CFa(\Surgnm)] \oplus [\{a_4\}\times \CFa(\Surgnm)] \oplus [\{y\}\times \CFKa(\companion)]$$
$$\CFKa(\WDnm,-1)\sim \{a_3\}\times \CFa(\Surgnm)$$

\noindent We use the symbol $\sim$ to mean that we have a bijection between the generators of e.g.  $\CFKa(\WDnm,1)$ and $ \{a_1\}\times \CFa(\Surgnm)$. It remains to understand the boundary operator for the chain complexes.  We direct our attention to the top filtration level for \WD.  We claim that, while {\em a priori} the chain complex $\CFKa(\WDnm,1)$ looks like $\CFa(\Surgnm)$ as stated above, it is in fact chain homotopy equivalent to the chain complex 
$$\CFMnm.$$
\noindent In other words, we have

\begin{theorem}
\label{thm:Meridian}

Let $\companion$ be a knot in $S^3$.  Then
$$\HFDnm \cong \HFMnm,$$
 \noindent Where the latter summation is taken over Spin$^c$ structures on $\Surgnm \cong \Z/t\Z$ and over filtration levels induced by $(\Surgnm,\mu_\companion).$

\end{theorem}

\noindent {\bf Remark:} The groups \HFMnospinc \ have a well-defined absolute $\Q$-grading which is a lift of a relative $\Z$-grading.  The isomorphism stated above does not, in general, hold on the level of graded abelian groups.  We will obtain the graded statement in the next section.

\begin{proof}
As previously noted, the generators of \CFD \ are of the form $\{a_1\}\times \CFa(\Surgnm)$.  By definition of the knot filtration, this implies that we can connect  $\{a_1\}\times \p$ to $\{a_1\}\times \q$ with a Whitney disk, $\phi$ for any $\p,\q\in \CFa(\Surgnm)$.  If we additionally require that $n_z(\phi)=n_w(\phi)=0$, then $\phi$ will be unique.  Recall that the Heegaard diagram hd$((\Surgnm,\mu_\companion)$ which went into the construction of hd$(\WDnm)$ came equipped with two basepoints, $z',w'$.  As described in \cite{HolDisk}, the point $w'$ induces a map
$$\spinc_{w'}  : \CFa(\Surgnm)\rightarrow  \mathrm{Spin}^c(\Surgnm),$$

\noindent which assigns to each generator in the chain complex a Spin$^c$ structure.  We first claim that if $\spinc_{w'}(\p)\ne \spinc_{w'}(\q)$ then the Whitney disk connecting  $\{a_1\}\times \p$ to $\{a_1\}\times \q$ will not admit any pseudo-holomorphic representatives.  This follows from the discussion in the proof of Lemma \ref{lemma:Lemma2}:  since $\spinc_{w'}(\p)\ne \spinc_{w'}(\q)$, the domain of the Whitney disk connecting $\{a_1\}\times \p$ to $\{a_1\}\times \q$ must contain the meridian $\mu_\companion$ in its boundary a non-zero number of times.  Since there are no corner points for the Whitney disk on the Heegaard diagram for the pattern (the disk has only degenerate corners at $\{x,a_1\}$ on hd$(Hopf)$), the domain of $\phi$ restricted to hd$(Hopf)$ is simply $n\cm \PerDom$, $n\ne0$ (where, as above $\PerDom$ indicates the periodic domain for hd$(Hopf)$).  Since $\PerDom$ has both positive and negative multiplicities (Figure \ref{fig:PerDom}) the disk cannot admit any pseudo-holomorphic representatives.

It follows that the boundary operator on \CFD \ respects the splitting along Spin$^c$ structures that it inherits as a set from $\CFa(\Surgnm)$.  It remains to understand the boundary operator for each Spin$^c$ structure.  Under the bijection between generators of \CFD \ and $\CFa(\Surgnm)$, we claim that boundary operator on \CFD \ is precisely the operator obtained from $\CFa(\Surgnm)$ by requiring $n_{z'}(\phi)=0$, in addition to $n_{w'}(\phi)=0$.  Since the Heegaard diagram for \Surg \ with both basepoints $z',w'$ is a compatible diagram for the knot $(\Surgnm,\mu_\companion)$, the theorem will follow from our claim and the definition of  \HFM.  

We prove the claim by examining the unique Whitney disk satisfying $n_z({\phi})=n_w({\phi})=0$ which connects $\{a_1\}\times \p$ to $\{a_1\}\times \q$ for any $\p,\q$ with $\spinc_{w'}(\p)=\spinc_{w'}(\q)$.   Since the disk has no corner points on hd$(Hopf)$, it restricts to $n\cm \PerDom$ on hd$(Hopf)$.   However, in order for $\phi$ to admit a holomorphic representative, $n=0$ since $\PerDom$ has positive and negative multiplicities.  Thus, the multiplicities of $\phi$, like those of $\PerDom$, are zero in a neighborhood of the region where we formed the connected sum in our construction of hd(\WD). In particular, $n_{z'}(\phi)=n_{w'}(\phi)=0$.  Conversely, any holomorphic disk connecting $\p$ to $\q$ in $\CFa(\Surgnm,\spinc_{w'}(\p))$ which satisfies $n_{z'}(\phi)=n_{w'}(\phi)=0$ can be extended to a disk connecting $\{a_1\}\times \p$ to $\{a_1\}\times \q$.  Thus the holomorphic disks that connect $\{a_1\}\times \p$ to $\{a_1\}\times \q$ for $\p,\q$ with $\spinc_{w'}(\p)=\spinc_{w'}(\q)$ correspond to holomorphic disks in $\CFa(\Surgnm,\spinc_{w'}(\p))$ with $n_{z'}(\phi)=n_{w'}(\phi)=0$ as we claimed.   (Here we are implicitly using the gluing theorem for pseudo-holomorphic disks used to prove stabilization invariance of Heegaard Floer homology, Section $11$ of \cite{HolDisk}.) 

\end{proof}

\section{Computation of \HFD \ for $|t|\gg 0$}
\label{sec:largen}
In this section we compute \HFD \ for all sufficiently large values of the twisting parameter, $t$.  We do this by showing that $\bigoplus \HFKa(\Surgnm,\mu_\companion,\spinc_i)$ - which was identified with \HFD \ in the previous section - is determined by the filtered chain homotopy type of $\CFKinf(\companion)$ in much the same way that $\HFa(\Surgnm,\spinc_i)$ is. (Recall that \ons \ present a formula for $\HFa(\Surgnm,\spinc_i)$ in terms of the filtered chain homotopy type of  $\CFKinf(\companion)$.  This is summarized below but see \cite{Knots} and also \cite{IntegerSurgeries,RationalSurgeries} for further details.)  Using essentially the method of \ons, we will prove the following theorem:

\begin{theorem}

	\label{thm:merid}	
Let $\companion\subset S^3$ be a knot.  Then for an appropriate labeling of Spin$^c$ structures, $\spinc_m$, $\exists \ T\gg0$ such that for all $t>T$ the following holds:

$$\HFKa_*(\Surgnm,\mu_\companion,\spinc_m) \cong H_{*+d_{_-}(m)}(\Filt(\companion,m)) \oplus H_{*-2m+d_{_-}(m)}(\Filt(\companion,-m-1))$$

$$\HFKa_*(S^3_{-t}(\companion),\mu_\companion,\spinc_m) \cong H_{*-d_{_+}(m)}(\frac{\CFa(S^3)}{\Filt(\companion,m)}) \oplus H_{*-2m-d_{_+}(m)}(\frac{\CFa(S^3)}{\Filt(\companion,-m-1)})$$

where 
$$
d_{\pm}(m)=\left(\frac{t-(2m\pm t)^2}{4t}\right).
$$ 
\end{theorem}

\noindent {\bf Remarks:} Here, as usual, the labeling of Spin$^c$ structures is determined by the condition that $\spinc_m$ can be extended over the cobordisms $-W_{t}'$ (resp. $W_{-t}$)  to yield a Spin$^c$ structure $\spincX$ satisfying:
$$ <c_1(\spincX_m),[S]> + t = 2m \ \ (resp. <c_1(\spincX_m),[S]> - t = 2m.) $$

\noindent $W_{t}$ denotes the cobordism from $S^3$ to $S^3_t(K)$ associated to the two-handle addition along $K$ with framing $t$.  A negative sign on a cobordism ($-W_{t}$) denotes the same cobordism with orientation reversed whereas a prime ($W_{t}'$) indicates that we ``turn the cobordism around'', viewing it as a cobordism from $-S^3_t(K)$ to $-S^3$. Thus $-W_{t}'$ is a cobordism from $S^3_t(K)$ to $S^3$.

\begin{proof} 
The theorem follows from an examination of the proof of Theorems $4.1$ and $4.4$ of \cite{Knots}. Recall that these theorems yield isomorphisms: \begin{eqnarray}
	\label{eq:iso}
\HFa_*(\Surgnm,\spinc_m) \cong H_{*+d_-(m)}(C\{\mathrm{max}(i,j-m)=0\}) \\
\HFa_*(S^3_{-t}(\companion),\mu_\companion,\spinc_m) \cong H_{*-d_+(m)}(C\{\mathrm{min}(i,j-m)=0\}).
\end{eqnarray}

\noindent Where $C\{\mathrm{max}(i,j-m)=0\}$ denotes the subquotient complex of the $\Z\oplus\Z$ filtered chain complex, $\CFKinf(\companion)$, generated by triples $[\x,i,j]$ with  $i$ and $j$ filtration indices satisfying the specified constraint.   The first isomorphism is induced by a chain map:
$$\Phi_{\spincX_m}\colon \CFa(S^3_{t},\spinc_m) \longrightarrow
C\{\mathrm{max}(i,j-m)=0\}$$

\noindent defined by counting pseudo-holomorphic triangles with appropriate boundary conditions in the $g$-fold symmetric product of $\Sigma_g$. The boundary conditions are determined by a doubly-pointed Heegaard triple diagram \Triple \ specifying the $2$-handle cobordism $-W_{t}'$.  The three three-manifolds specified by the triple diagram are \Yab$=$\Surg, \Yag$=S^3$, \Ybg$=\#^{g-1}S^1\times S^2$ (see Figure \ref{fig:Surgery}).   The map is defined by:
$$\Phi_{\spincX_m}[\x]=\sum_{\y\in\Ta\cap\Tc}
\sum_{\{\psi\in\pi_2(\x,\Theta,\y)\big| 
n_\FiltPt(\psi)-n_\BasePt(\psi)=m-\Filt(\y),
\Mas(\psi)=0\}}
\left(\#\ModFlow(\psi)\right)\cm [\y,-n_\BasePt(\psi),m-n_\FiltPt(\psi)],$$ 

\noindent where $\Theta$ is a top-dimensional generator for $\HFa(\#^{g-1}S^1\times S^2)$.
The second isomorphism is induced by a similar map which goes in the opposite direction. The condition on the homotopy classes of triangles, $\psi$, in the above map ensures that the image of the map is $C\{\mathrm{max}(i,j-m)=0\}$ (this last statement follows from Equation \eqref{eq:SpinCFormula} found in the proof of Theorem \ref{thm:largen} below).

We would like to refine \ons's theorem to determine the knot Floer homology of $(\Surgnm,\mu_\companion)$.  In fact, there is a natural $2$-stage sequence of subcomplexes:
$$0\subseteq C\{i<0,j=m\}\subseteq C\{\mathrm{max}(i,j-m)=0\}.$$
\noindent Viewing this sequence as a $2$-step filtration, we claim that its filtered chain homotopy type is filtered chain homotopy equivalent to that of the filtration of $\CFa(\Surgnm,\spinc_m)$ induced by $(\Surgnm,\mu_\companion)$.  Summing over Spin$^c$ structures on \Surg, Theorem \ref{thm:merid} will follow immediately, since
$$C_*\{i<0,j=m\}\simeq \Filt_{*-2m}(K,-m-1)$$  
$$\frac{C_*\{\mathrm{max}(i,j-m)=0\}}{ C_*\{i<0,j=m\}} \simeq C_*\{i=0,j\le m\} \simeq \Filt_*(K,m)$$
\noindent by the discussion in Section $3.5$ of \cite{Knots}. (here $\simeq$ indicates chain homotopy equivalence). The case when $t<0$ is similar.

In proving our claim, the key observation will be that the triple diagram \Triple \ used to define $\Phi_{\spincX_m}$  not only specifies a Heegaard diagram for the knot $(S^3,\companion)$ (as was used in the proof of \ons's Theorems) but also a Heegaard diagram for the knot $(\Surgnm,\mu_\companion)$, with the addition of a basepoint, $z'$.  Indeed, $(\Sigma,\alphas,\betas,w,z')$ is a diagram for $(\Surgnm,\mu_\companion)$ as was noted in the discussion of the Heegaard diagram for the companion in Section \ref{sec:heegs}.  We show the local picture of the Heegaard triple near the basepoints in Figure \ref{fig:Surgery}.
\begin{figure}
	\begin{center}
\psfrag{b}{$\beta_g$}
\psfrag{a}{$\alpha_g$}
\psfrag{c}{$\gamma_g$}
\psfrag{theta}{$\Theta$}
\psfrag{w}{$w$}
\psfrag{z}{$z$}
\psfrag{y}{$y$}
\psfrag{z'}{$z'$}
\psfrag{x-1}{$x_1$}
\psfrag{x-2}{$x_2$}
\psfrag{x-3}{$x_3$}
\psfrag{x0}{$x_0$}
\psfrag{x1}{$x_{-1}$}
\psfrag{x2}{$x_{-2}$}
\psfrag{x3}{$x_{-3}$}
\psfrag{0}{$0$}
\psfrag{1}{$1$}
\psfrag{2}{$2$}
\psfrag{psi1}{$\psi_1$}
\psfrag{psi2}{$\psi_{-2}$}
\psfrag{alpha1}{$\alpha_1$}
\caption{\label{fig:Surgery}
Local picture of Heegaard triple for $-W_t'$ near the basepoints (here $t=6$).  We refer to this part of the Heegaard diagram as the ``winding region''. In blue shade is the domain of a small triangle $\psi_1$ connecting $x_1$ to $x_0$.  In red shade is a small triangle $\psi_{-2}$ connecting $x_{-2}$ to $x_0$, with $n_w(\psi_{-2})=2$. Note that while the picture looks like the diagram for hd$(\WDnm)$ near the connect sum region,  this is a Heegaard triple diagram, and we make no reference (for the moment) to the diagram for the Whitehead doubles.  Thus the notation here is independent of the previous two sections.} 
\includegraphics{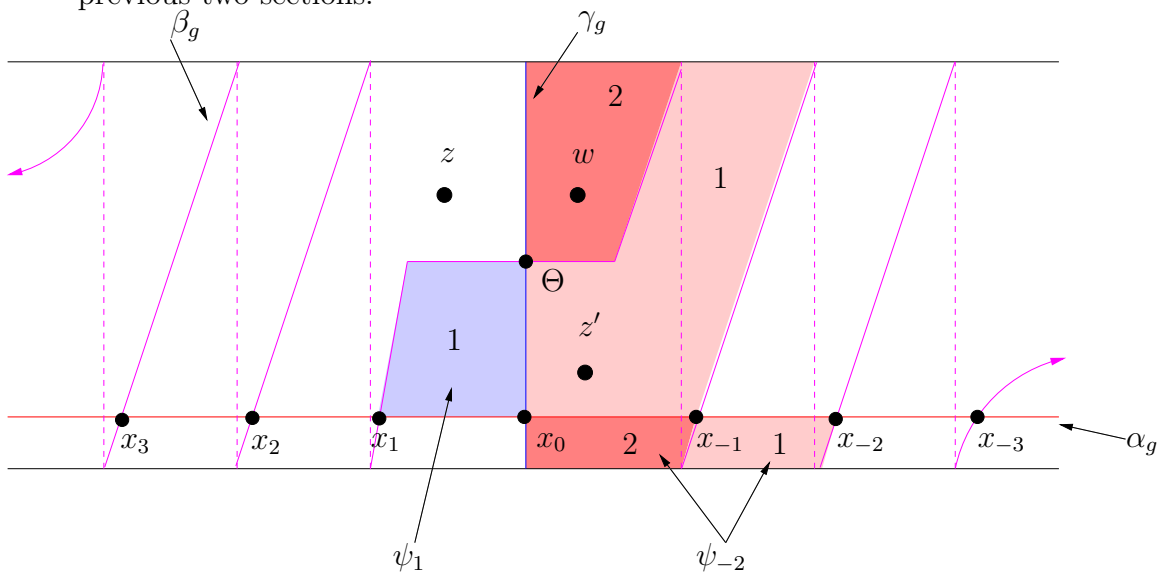}
\end{center}
\end{figure}

An intersection point $\x\in \Ta\cap\Tb$ is said to be supported in the winding region if the component of $\x$ in $\alpha_g$ lies in the local picture of Figure \ref{fig:Surgery}.  Intersection points in the winding region are in $t$ to $1$ correspondence with intersection points $\y \in \Ta\cap\Tg$.  When $t$ is sufficiently large, the pigeonhole principle shows that there exist an entire Spin$^c$ equivalence class of intersection points for $\CFa(\Surgnm)$ supported in the winding region.  Furthermore, the generators of this equivalence class are in bijection with generators of \newline $C\{\mathrm{max}(i,j-m)=0\}$.  As shown by \ons, this bijection is induced by canonical ``small'' triangles supported entirely in the winding region which connect generators of the two complexes, see Figure \ref{fig:Surgery}.  In fact the pseudo-holomorphic representatives of these triangles constitute the highest order terms of $\Phi_{\spincX_m}$ (with respect to the filtration given by negative area of triangles), and subsequently induce the isomorphism on homology given by Equation \eqref{eq:iso}.
 
From Figure \ref{fig:Surgery} we see that there is a unique intersection point $x_0\in \alpha_g\cap \gamma_g$. The multiplicity of each small triangle at the basepoints measures how far to the right or left of $x_0$ the $\alpha_g$ component of the corresponding generator of $\CFKa(\Surgnm,\mu_\companion,\spinc_m)$ lies.  Suppose $\psi\in \pi_2(\x,\Theta,\y)$ is a small triangle. If $k=n_z(\psi)\ge0$ then the $\alpha_g$ component of $\x$ is $x_k$, and  $[\x]$ is mapped by $\Phi_{\spincX_m}$ to the quotient complex, $C\{0,j\le m\}$, of $C\{\mathrm{max}(i,j-m)=0\}$.  If, on the other hand, $l=n_w(\psi)>0$ the $\alpha_g$ component of $\x$ is $x_{-l}$ and  $[\x]$ is mapped  to the subcomplex, $C\{i<0,j=m\}$.

We now claim that the filtration difference between two generators of $\CFKa(\Surgnm,\mu_\companion,\spinc_m)$ is equal to $\pm1$ if their $\alpha_g$ components lie on opposite sides of $x_0$, and is zero otherwise.  To see this, observe that the boundary of the domain of the unique Whitney disk, $\phi$, connecting  $\x,\y\in \CFKa(\Surgnm,\mu_\companion,\spinc_m)$ with $n_w(\phi)=0$ contains the arc on $\beta_g$ connecting $x_1$ and $x_{-1}$ if, and only if, the $\alpha_g$ components of $\x$ and $\y$ are on opposite sides of $x_0$.  Furthermore, the multiplicity with which this arc occurs in $\partial\phi$ is $1$ if $\x$ and $\y$ are on opposite sides of $x_0$ and $0$ otherwise.  This proves the claim.

Our original claim about the identification of the $2$-stage filtration of \newline $C\{\mathrm{max}(i,j-m)=0\}$ with the filtration induced by $(\Surgnm,\mu_\companion)$ follows from the preceding two paragraphs. Explicitly, the first paragraph showed that under the bijection between generators of $\CFKa(\Surgnm,\mu_\companion,\spinc_m)$ and $C\{\mathrm{max}(i,j-m)=0\}$ induced by the small triangles, the generators of $\CFKa(\Surgnm,\mu_\companion,\spinc_m)$ with $\alpha_g$ component lying to the left of $x_0$ correspond exactly to $C\{i=0,j\le m\}$ while the  generators with $\alpha_g$ component to the right of $x_0$ correspond to $C\{i<0,j=m\}$. The second paragraph then identified the generators with $\alpha_g$ component to the left of $x_0$ with the top filtration level of $\CFKa(\Surgnm,\mu_\companion,\spinc_m)$ and those with $\alpha_g$ component to the right of $x_0$ with the bottom filtration level.   This completes our proof of Theorem \ref{thm:merid}.

\end{proof}

\subsection{Gradings}
\label{subsec:gradings}
We conclude this section by determining the absolute Maslov grading of \HFD, when $|t|\gg0$.  We show the following:

\begin{theorem} 
	\label{thm:largen}	With notation as above, there are isomorphisms of absolutely $\Z$-graded abelian groups for all $t>T$:
	
	$$\HFDgnm \cong  \bigoplus_{m=\lfloor -\frac{t}{2} \rfloor +1}^{m=\lfloor \frac{t}{2} \rfloor }[ H_{*-1}(\Filt(\companion,m)) \oplus H_{*-1}(\Filt(\companion,-m-1))]$$

$$\HFDgnmneg \cong  \bigoplus_{m=\lfloor -\frac{t}{2} \rfloor +1}^{m=\lfloor \frac{t}{2} \rfloor }[ H_{*}(\frac{\CFa(S^3)}{\Filt(\companion,m)}) \oplus H_{*}(\frac{\CFa(S^3)}{\Filt(\companion,-m-1)})].$$	
\end{theorem}

\begin{proof}
	By a sequence of Heegaard moves, each of which occur in the complement of the basepoint $w$, we can convert 
	hd$(\WDnm)$ to hd$(\companion)$. See Figure \ref{fig:Handleslide}.  Since these moves occur in the complement of $w$, the absolute grading of any generators unaffected by the Heegaard moves is unchanged throughout the process.  It follows that intersection points of the form:
	$$\{y\} \times \p_i,  \in \CFKa(P) \times \CFKa(\companion)$$
	\noindent inherit the grading which $\p_i$ has, thought of as a generator of $\CFKa(\companion)$.  As mentioned in the proof of Theorem \ref{thm:merid} above, intersection points $\p_i\in \CFKa(\companion)$ are in $1$ to $t$ correspondence with intersection points in the winding region, and hence with intersection points of the form:
	$$\{a_j\} \times \q_i \in \{a_j\}\times\CFa(\Surgnm),$$
	\noindent where $q_i$ is supported in the winding region.  More explicitly, each point of $\{y\}\times \CFKa(\companion)$ is of the form $\{y\}\times\{x_0,\p\}$, for some $(g-1)$-tuple of intersection points, $\p$, while each point of $\{a_j\}\times \CFa(\Surgnm)$ is of the form $\{a_j\}\times \{x_k,\p\}$ for some $x_k \in\beta_g\cap\lambda_\companion\#\lambda_P$, where $k=\lfloor -\frac{t}{2} \rfloor +1,\ldots,\lfloor \frac{t}{2} \rfloor$.  In order to determine the absolute gradings for the Floer homology of the Whitehead double, we first understand the grading on the intersection points supported in the winding region. These points are partitioned into two groups - those points whose $\beta_g$ component is to the left of $x_0$, and those whose $\beta_g$ component is to the right.  We first handle those points to the left of $x_0$.

	\begin{lemma} 
		\label{lemma:rectangle}
		Let $k>0$.  Then for $t\ll0$ we have
	$$gr(\{a_1\}\times \{x_k,\p\}) = gr(\{y\}\times \{x_0,\p\})+1,$$
	\noindent while for $t\gg0$,
	$$gr(\{a_4\}\times \{x_k,\p\}) = gr(\{y\}\times \{x_0,\p\})-1.$$
\end{lemma}
\begin{proof} For $k>0$, we can complete the small triangles connecting $x_k$ to $x_0$ used in the proof of Theorem \ref{thm:merid} with domains on the diagram for $P$ to obtain domains of Whitney disks.  This is shown in Figure \ref{fig:positiverectangle} for $t\gg 0$.  The Whitney disks connect $\{a_1\}\times \{x_k,\p\}$ to $\{y\}\times \{x_0,\p\}$ if $t\ll0$ and $\{y\}\times \{x_0,\p\}$ to $\{a_4\}\times \{x_k,\p\}$ if $t\gg0$.  The Maslov index is easily calculated (using, for example Corollary $4.10$ of \cite{Lipshitz}) to be $1$ in both cases. 
\end{proof}

\begin{figure}
	\begin{center}
\psfrag{b1}{$\beta_1$}
\psfrag{a}{$\alpha_g$}
\psfrag{c}{$\gamma_g$}
\psfrag{theta}{$\Theta$}
\psfrag{w}{$w$}
\psfrag{z}{$z$}
\psfrag{K}{\tiny{N}}
\psfrag{K-1}{\tiny{N-1}}
\psfrag{K-2}{\tiny{N-2}}
\psfrag{K-3}{\tiny{N-3}}
\psfrag{-K+1}{\tiny{-N+1}}
\psfrag{y}{$y$}
\psfrag{a4}{$a_4$}
\psfrag{x-1}{$x_1$}
\psfrag{x0}{$x_0$}
\psfrag{mk}{$\mu_\companion\#\mu_P$}
\psfrag{l}{$\lambda_\companion\#\lambda_P$}
\psfrag{2}{$2$}
\psfrag{psi1}{$\psi_1$}
\psfrag{psi2}{$\psi_{-2}$}
\psfrag{alpha1}{$\alpha_1$}
\psfrag{Front}{Front}
\psfrag{Back}{Back}
\psfrag{C}{Companion}
\psfrag{Pf}{Pattern Front}
\psfrag{Pb}{Pattern Back}
\psfrag{m}{$\mu$}
\caption{\label{fig:positiverectangle}
Illustration for Lemmas \ref{lemma:filt} and \ref{lemma:rectangle}.  The Figure shows how to complete a Whitney triangle $\psi\in \pi_2(\p,\Theta,\{x_0,\q_i\})$ with appropriate multiplicities near $\Theta$ on the Heegaard triple corresponding to $-W_t'$  to a Whitney disk $\phi\in \pi_2(\{y\}\times \{x_0,\q_i\},\{a_4\}\times \p)$ on the Heegaard diagram hd$(\WDnm)$.  If the triangle has Maslov index $k$, the corresponding Whitney disk will have index $k+1$.  Multiplicities of the domain of $\phi$ are shown in the shaded regions. When the framing $t$ is negative, a similar procedure can complete a triangle to a disk $\phi\in \pi_2(\{a_1\}\times \p,\{y\}\times \{x_0,\q_i\})$  } 
\includegraphics{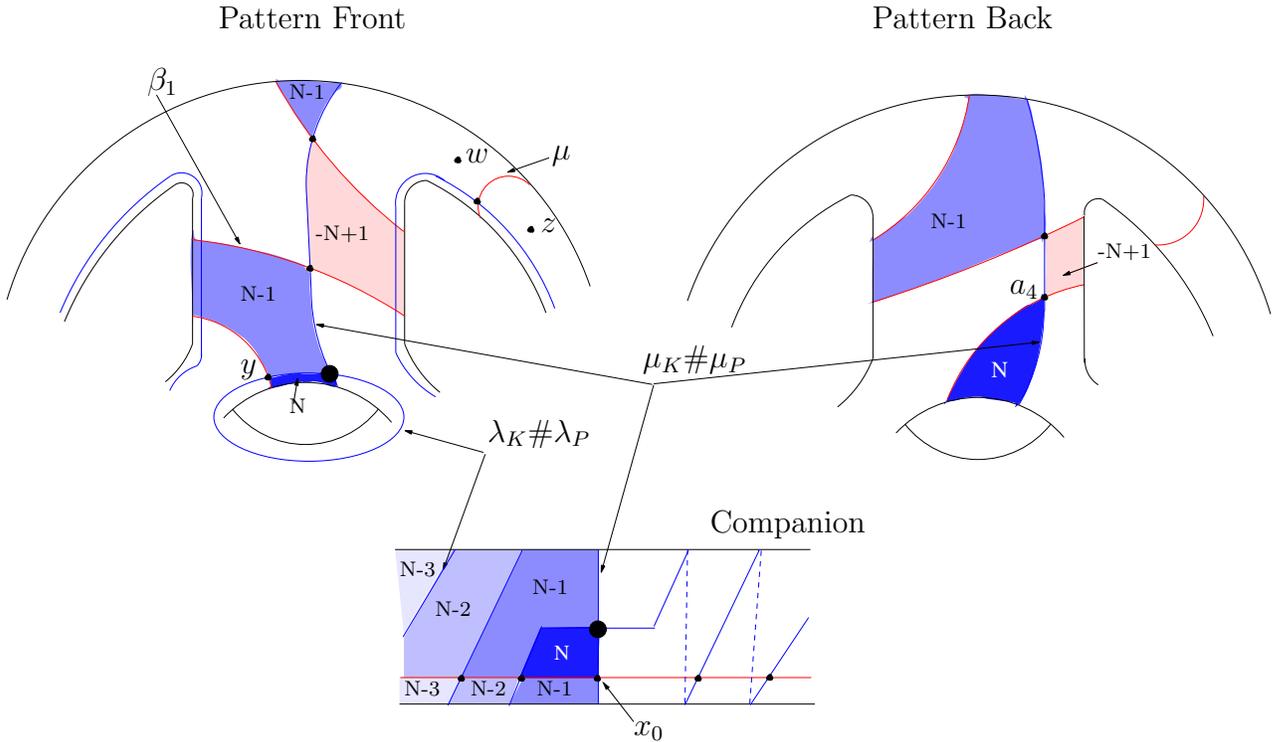}
\end{center}
\end{figure}

Next we deal with the points to the right of $x_0$.   Before stating the lemma, we remark that each generator of $\CFKa(K)$ is of the form $x_0\times\p$, and hence we can think of the $(g-1)$-tuple $\p$ as having a filtration, $\Filt(\p)$ inherited from the filtration of $\CFKa(\companion)$.  Note also that each point in the winding region $x_k\times \p$ uniquely corresponds to an intersection point $x_0\times \p\in \CFKa(\companion)$.
\begin{lemma}
	\label{lemma:rectangle2}
	Let $-k<0$. Then for $t\ll0$ we have
	$$gr(\{a_1\}\times \{x_{-k},\p\}) = gr(\{y\}\times \{x_0,\p\})+2\Filt(\p)+1,$$
	\noindent while for $t\gg0$,
	$$gr(\{a_4\}\times \{x_{-k},\p\}) = gr(\{y\}\times \{x_0,\p\})+2\Filt(\p)-1.$$
\end{lemma}
\begin{proof} We handle only the second case, as the first is similar. We would like to complete the small triangle $\psi_{-k}$ connecting $\{x_k,\p\}$ to $\{x_0,\p\}$ by a domain on hd$(Hopf)$, as in the previous lemma.  However, since the small triangle is now supported to the right of $x_0$, the multiplicities of the domain of the triangle near the connect sum region (the black hole) are not suitable for completion.  Thus we pick a homotopy class of triangles connecting $\{x_k,\p\}$ to $\{x_0,\p\}$ which has multiplicity $0$ in the two domains on the right of the connect sum tube.  Since $\psi_{-k}$ has multiplicity $k$ at the basepoint $w$, we first subtract off $k\cm [\Sigma]$ from the domain of $\psi_{-k}$.  This has the effect of lowering the Maslov index of $\psi_{-k}$ by $2k$.  Next, we subtract the generator of the space of triply-periodic domains shown in Figure \ref{fig:triply}, which we denote by $\PerDom^3$.  The effect that this has on the Maslov index of the triangle can be determined from the definition of the absolute grading of Floer homology for torsion Spin$^c$ structures (Equation $(12)$ of \cite{HolDiskFour}:
	\begin{equation}
		\label{eq:Maslov}
	\mu(\psi_{-k}-\PerDom^3)-\mu(\psi_{-k})=\frac{1}{4}( (c_1(\spinc_w(\psi_{-k}-\PerDom^3))^2-(c_1(\spinc_w(\psi_{-k}))^2).
\end{equation}
	\noindent In the above equation, $\spinc_w(\psi)$ denotes the Spin$^c$ structure on the two handle cobordism $-W_t'$ associated to the triangle $\psi$ via the basepoint $w$, and $c_1$ denotes its first Chern class (for a description of how a choice of basepoint and homotopy class of triangles specifies a Spin$^c$ structure on $-W_t'$, see Section $8$ of \cite{HolDisk}).
	Now for an arbitrary homotopy class of triangles $\psi$ connecting $ \{x_k,\p\}$ to $\{x_0,\p\}$, an analogue of  Equation $(14)$ of \cite{Knots} states:
\begin{equation}
\label{eq:SpinCFormula}
\langle c_1(\spincrel(\{x_0,\p\})),[{\widehat F}]\rangle +2(n_{\FiltPt}(\psi)-n_{\BasePt}(\psi))-t
= \langle c_1(\spinc_{\BasePt}(\psi)),[S]\rangle.
\end{equation}
\noindent Here $\spincrel(\{x_0,\p\}$ denotes the Spin$^c$ structure on the zero surgered manifold, $S^3_0(K)$, associated to $\{x_0,\p\}$ by the basepoint $w$.  $[{\widehat F}]$ denotes the homology class in $H_2(S^3_0(K),\Z)$ corresponding to the Seifert surface of $K$, capped off by the meridian disk of the surgery torus, and $[S]$ denotes the generator of $H_2(-W_t',\Z)$.  Recall that for knots in $S^3$, the knot Floer homology filtration, $\Filt$, can be thought of as a filtration by Spin$^c$ structures on $S^3_0(K)$.  Under this correspondence, 
$$\langle c_1(\spincrel(\{x_0,\p\})),[{\widehat F}]\rangle = 2\Filt(\{x_0,\p\}):= 2\Filt(\p),$$
\noindent and hence for the small triangle $\psi_{-k}$ Equation \eqref{eq:SpinCFormula} becomes:
$$ 2\Filt(\p)-2k-t = \langle c_1(\spinc_{\BasePt}(\psi_{-k})),[S]\rangle. $$
By subtracting $\PerDom^3$ from $\psi_{-k}$ we change the Spin$^c$ structure associated to the triangle:
$$\langle c_1(\spinc_{\BasePt}(\psi_{-k}-\PerDom^3)),[S]\rangle = \langle c_1(\spinc_{\BasePt}(\psi_{-k}))-2\mathrm{PD}[S],[S]\rangle=\langle c_1(\spinc_{\BasePt}(\psi_{-k})),[S]\rangle+2t$$

We can now compute $c_1^2$ for the Spin$^c$ structures associated to $\psi_{-k}$ and $\psi_{-k}-\PerDom^3$ and determine the difference in their Maslov indices using Equation \eqref{eq:Maslov}:
$$	\mu(\psi_{-k}-\PerDom^3)-\mu(\psi_{-k})=\frac{1}{4}( \frac{(2\Filt(\p)-2k+t)^2}{-t}-\frac{(2\Filt(\p)-2k-t)^2}{-t})= 2(k-\Filt(\p))$$

Thus, we arrive at the triangle $\psi_{-k}'=\psi_{-k}-k\cm[\Sigma]-\PerDom^3$ connecting $\{\x_k,\p\}$ to $\{x_0,\p\}$ whose multiplicities near the connect sum region are shown in Figure \ref{fig:triply}.  The Maslov index of this triangle is:
$$\mu(\psi_{-k}-k\cm[\Sigma]-\PerDom^3)=\mu(\psi_{-k})-k\cm\mu([\Sigma])-\mu(\PerDom^3)=0-2k+2(k-\Filt(\p))=-2\Filt(\p).$$

The domain of $\psi_{-k}'$ can be completed on the Heegaard diagram for the pattern to yield the domain of a Whitney disk $\phi\in \pi_2(\{y\}\times \{x_0,\p\}, \{a_4\}\times \{x_k,\p\})$.  As in Lemma  \ref{lemma:rectangle} above, the Maslov index of this disk is one higher than the Maslov index of $\psi_{-k}'$. That is, 
$$\mu(\phi)=\mu(\psi_{-k})+1=-2\Filt(\p)+1.$$  
Since the relative grading in knot Floer homology is determined by $gr(\x)-gr(\y)=\mu(\phi)-2n_w(\phi)$ for $\phi\in\pi_2(\x,\y)$, the proof of the lemma is completed.
\end{proof}

\begin{figure}
	\begin{center}
\psfrag{b}{$\beta_g$}
\psfrag{a}{$\alpha_g$}
\psfrag{c}{$\gamma_g$}
\psfrag{theta}{$\Theta$}
\psfrag{w}{$w$}
\psfrag{z}{$z$}
\psfrag{a4}{$a_4$}
\psfrag{x2}{$x_{-2}$}
\psfrag{x0}{$x_0$}
\psfrag{2}{$2$}
\psfrag{psi}{$\psi_{-2}$}
\psfrag{psi'}{$\psi_{-2}'=\psi_{-2}-2\cm[\Sigma_g]-\PerDom^3$}
\psfrag{0}{$0$}
\psfrag{1}{$1$}
\psfrag{2}{$2$}
\psfrag{3}{$3$}
\psfrag{-4}{$-4$}
\psfrag{5}{$5$}
\psfrag{-1}{$-1$}
\psfrag{-2}{$-2$}
\psfrag{-3}{$-3$}
\psfrag{T}{$t$}
\psfrag{T-1}{$t-1$}
\psfrag{T-2}{$t-2$}
\psfrag{T-3}{$t-3$}
\psfrag{T-4}{$t-4$}
\psfrag{T-5}{$t-5$}
\psfrag{-T}{$-t$}
\psfrag{-T+1}{$-t+1$}
\psfrag{-T+2}{$-t+2$}
\psfrag{-T+6}{$-t+6$}
\psfrag{-T+3}{$-t+3$}
\psfrag{-T+4}{$-t+4$}
\psfrag{-T+5}{$-t+5$}
\psfrag{-T-1}{$-t-1$}
\psfrag{P}{$\PerDom^3$}
\caption{\label{fig:triply}  
Illustration for Lemma \ref{lemma:rectangle2}.  The first part of the Figure depicts multiplicities of the small triangle $\psi_{-2}$.  The second shows multiplicities of the generator of the space of triply-periodic domains, $\PerDom^3$. The last part shows multiplicities of the triangle obtained from $\psi_{-2}$ by subtracting $2\cm[\Sigma]$ and $\PerDom^3$. The domain of $\psi_{-2}'$ can be completed to the domain of a Whitney disk on the Heegaard diagram for the Whitehead double, as in the proof of Lemma \ref{lemma:rectangle}.   } 
\includegraphics{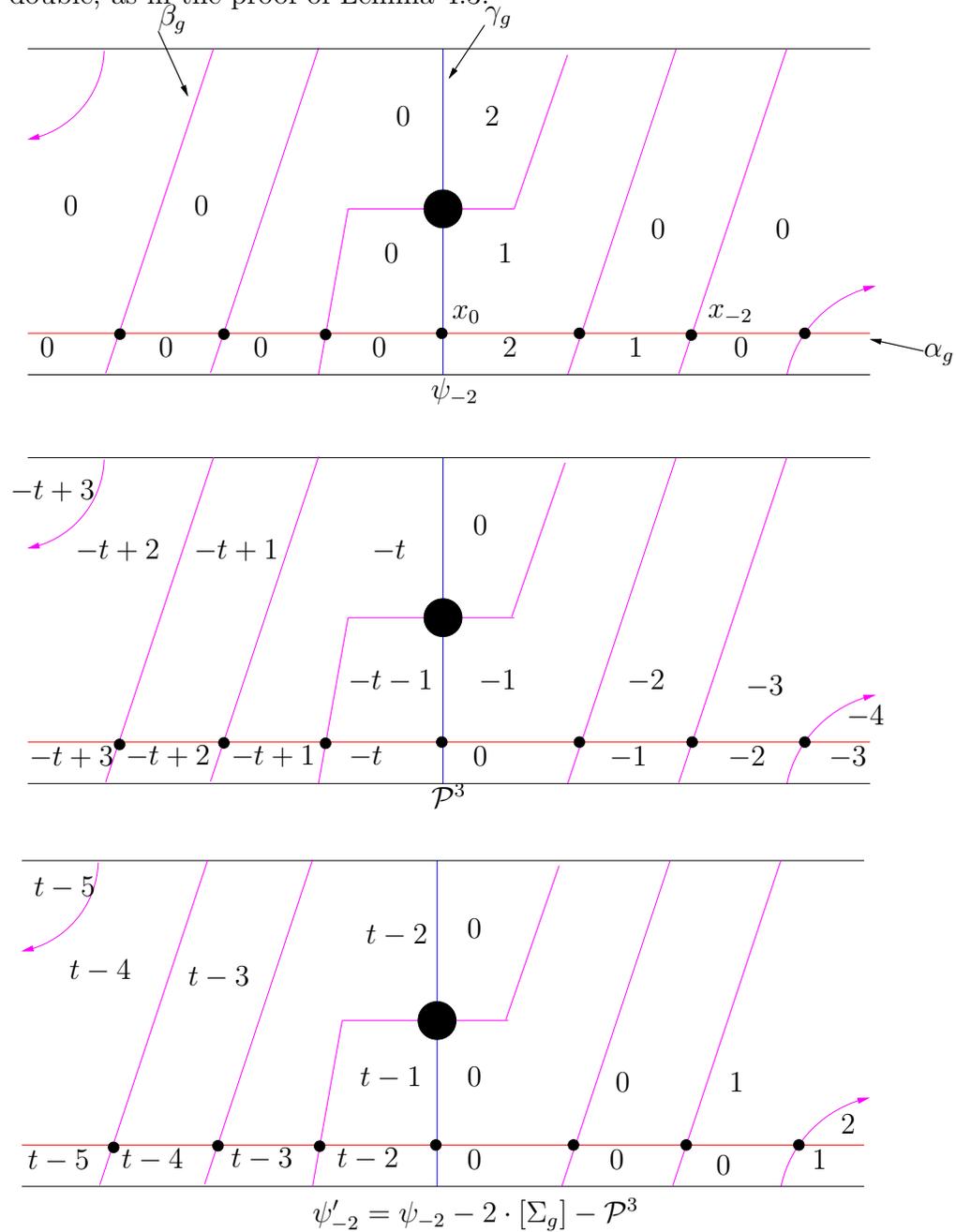}
\end{center}
\end{figure}

For intersection points generating \HFD \ which are supported in the winding region, the two lemmas above are enough to complete the proof of Theorem \ref{thm:largen}.  More explicitly, the isomorphism in Theorem \ref{thm:merid},
$$\HFKa_*(\Surgnm,\mu_\companion,\spinc_m) \cong H_{*+d_{_-}(m)}(\Filt(\companion,m)) \oplus H_{*-2m+d_{_-}(m)}(\Filt(\companion,-m-1)),$$
was proved by looking at intersection points in the winding region, with the $\Filt(\companion,m)$ summand corresponding to those $\{x_k,\p\}\in \CFKa(\Surgnm)$ with $x_k$ to the left of $x_0$, and $\Filt(\p)\le m$.  By Lemma \ref{lemma:rectangle} the grading of these intersection points is the same as the grading of $\{x_0,\p\}$ (or shifted up by $1$, depending on whether $t$ is positive or negative) and hence the same as that of $\Filt(\companion,m)$.  The correction factor $d_{_-}(m)$ is missing since the Whitehead double is knot in $S^3$ and $d_{_-}(m)$ was inherited from the grading of \Surg.   On the other hand, the  $\Filt(\companion,-m-1)$ summand corresponded to those $\{x_{-l},\p\}\in \CFKa(\Surgnm)$ with $x_{-l}$ to the right of $x_0$, and $\Filt(\p)>m$, i.e.  $C\{i<0,j=m\}$.  Here we relied on the chain homotopy equivalence:
$$C_*\{i<0,j=m\}\simeq \Filt_{*-2m}(K,-m-1).$$  
\noindent This chain homotopy equivalence is the same as that used in the proof of Proposition $3.8$ of \cite{Knots}, which identifies the Floer homology of a knot $K$, and its reverse $-K$. Under this chain homotopy, the filtration of a generator $\Filt(\p)$ is sent to $-\Filt(\p)$, and the grading is shifted by $-2\Filt(\p)$.  This accounts for the grading shift of $2\Filt(\p)$ seen in Lemma \ref{lemma:rectangle2}, while the absence of the correction term $-2m$ is due to the fact that the basepoint $\BasePt$ used in the calculation of the grading of \Surg ~ is absent from the diagram for the Whitehead double.

We would like to claim that we are done with the calculation of the absolute grading for \HFD.  However, in the proof of Theorem \ref{thm:merid} we relied on the fact that we could calculate $\HFKa_*(\Surgnm,\mu_\companion,\spinc_m)$ for each Spin$^c$ structure separately, each time looking only at points in the winding region.  In the case of the Whitehead double all the Spin$^c$ structures on \Surg ~ are grouped into \HFD, and hence we must determine the absolute gradings of the generators of \HFD ~ lying outside the winding region.

\begin{figure}
	\begin{center}
\psfrag{b}{$\beta_g$}
\psfrag{a}{$\alpha_g$}
\psfrag{c}{$\gamma_g$}
\psfrag{w}{$w$}
\psfrag{0}{0}

\psfrag{M}{M}
\psfrag{N}{N}
\psfrag{-2N}{-2N}
\psfrag{-5N}{-5N}
\psfrag{N+M+1}{N+M+1}
\psfrag{5}{5}
\psfrag{M-N-1}{M-N-1}
\caption{\label{fig:badtriangle}
Illustration of the domain of an arbitrary triangle, $\psi$, connecting a generator supported away from the winding region to a generator of $\HFa(S^3)$.  Since $n_w(\psi)=0$, we see that in order for $\psi$ to be holomorphic, $N=0$.  Otherwise the domain of $\psi$ would have negative multiplicity.  With $N=0$ we can complete $\psi$ to a Whitney disk as in Lemmas \ref{lemma:rectangle} and \ref{lemma:rectangle2}. } 
\includegraphics{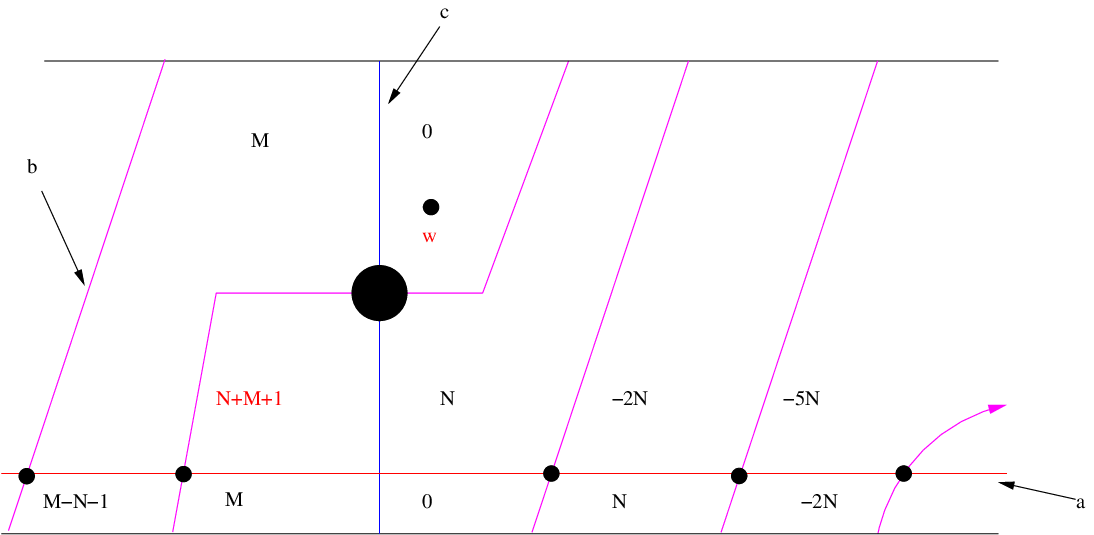}
\end{center}
\end{figure}

To handle these generators we argue as follows: by letting $|t|\gg 0$ be sufficiently large, we can ensure that every Spin$^c$ structure on \Surg ~ with $$\HFKa(\Surgnm,\mu_\companion,\spinc_m)\ne \Z$$ is generated by intersection points supported in the winding region.  In other words, we make the framing large enough so that the intersection points outside of the winding region only lie in Spin$^c$ structures on \Surg ~ for which the groups $\HFKa_*(\Surgnm,\mu_\companion,\spinc_m)$ have stabilized  (we know that these groups must eventually all be $\Z$, by the adjunction inequality).  For these Spin$^c$ structures, the map $\Phi_{\spincX_m}$ induces an isomorphism 
$$\HFKa_*(\Surgnm,\mu_\companion,\spinc_m)\cong \HFa(S^3)=\Z.$$
\noindent  Now since the maps $\Phi_{\spincX_m}$ are invariants of the cobordism $-W_{t}'$ and Spin$^c$ structure $\spincX_m$, it follows that we can calculate them with an arbitrary Heegaard triple diagram.  Since 
$\Phi_{\spincX_m}$ is an isomorphism for the $\spinc_m$ in the stable range, it follows that there exists a pseudo-holomorphic Whitney triangle $\psi$ with $\mu(\psi)=0$ connecting the generator of $\HFKa(\Surgnm,\mu_\companion,\spinc_m)$ to the generator of $\HFa(S^3)$.  This implies that the multiplicities of the domain of $\psi$ must all be positive and it follows that the multiplicities of $\psi$ in the domains to the right of the connect sum region must all be zero (otherwise there would be negative multiplicity somewhere in the winding region, see Figure \ref{fig:badtriangle}).  We can complete this $\psi$ to a Whitney disk $\phi$ with $\mu(\phi)=1$ as in Lemma \ref{lemma:rectangle} and thus the proof of Theorem \ref{thm:largen} is complete.
\end{proof}

\section{Computation of $\widehat{HFK}(D_+(K,t),1)$ and $\tau(D_+(K,t))$ for all $t$}
In this section we use the skein exact sequence for knot Floer homology to interpolate between the case when $t\gg0$ and the case when $t\ll0$.  This will enable us to determine \HFD \ for all values of $t$.  We will also determine $\tau(\WDnm)$.   The analysis of the skein sequence will be similar to the technique used in \cite{STau} in the special case of Whitehead doubles of the $(2,2n+1)$ torus knots.  The main result of this section will be:
\begin{prop}
	Let $K\subset S^3$ be a knot with Seifert genus $g(K)=g$.  Then for $t\ge 2\tau(K)$ we have:
$$\HFKa_*(\WDnm,1)\cong \Z_{(1)}^{t-2g-2} \bigoplus_{i=-g}^{g}[H_{*-1}(\Filt(K,i))]^2, $$
and $\tau(\WDnm)=0$.	For $t< 2\tau(K)$ the following holds:
$$\HFKa_*(\WDnm,1)\cong \Z_{(1)}^{2\tau(K)-2g-2} \oplus\Z_{(0)}^{2\tau(K)-t} \bigoplus_{i=-g}^{g}[H_{*-1}(\Filt(K,i))]^2, $$
and  $\tau(\WDnm)=1.$
\end{prop}
	\noindent {\bf Remark:} This takes care of the top and bottom group for the Whitehead double, by the symmetry of knot Floer homology about $\Filt=0$.  It also proves Theorem \ref{thm:tau} stated in the introduction.

\begin{proof} It will be helpful to first rephrase Theorem \ref{thm:largen}:
\begin{theorem} 
	\label{thm:largen2}
	Let $\companion\subset S^3$ be a knot, and suppose $g(\companion)=g$ denotes the Seifert genus of $\companion$. Then for all $t>T>0$ there are isomorphisms of absolutely $\Z$-graded abelian groups:
	$$\HFDgnm \cong \Z_{(1)}^{t-2g-2} \bigoplus_{i=-g} ^{g}[ H_{*-1}(\Filt(\companion,i))]^2$$
$$\HFDgnmneg \cong  \Z_{(0)}^{|t|-2g}\bigoplus_{i=-g} ^{g}[ H_{*}(\frac{\CFa(S^3)}{\Filt(\companion,i)})]^2$$	
\end{theorem}
\begin{proof} This follows from the adjunction inequality for knot Floer homology, which implies that $H_*(\Filt(\companion,i))\cong0$ for $i<g$ and $H_*(\Filt(\companion,i))\cong Z_{(0)}$ for $i\ge g$.
\end{proof}

As in \cite{STau}, we note that it is possible to change \WD \ to \WDneg \ by a sequence of $2t$ crossing changes, each of which change a negative crossing in the twisting region to a positive crossing.  Theorem $8.2$ of \cite{Knots} shows that corresponding to each crossing change, there is a long exact sequence relating the knot Floer homology groups of \WD, \WDm, and the two component link obtained from the oriented resolution of the crossing which we change.  In each case, this link is the positive Hopf link, which we denote by $H$. Summarizing the discussion of \cite{STau}, the skein exact sequence for the top filtration level takes the following form: $$
\begin{CD}
	... @>>>	\HFDnm
	@>{f_1}>> \F_{(\OneHalf)} @>{f_2}>> 	 \HFDnmm @>{f_3}>>...
\end{CD}
$$

\noindent Here the maps $f_1$ and $f_2$ lower homological degree by one-half and $f_3$ is non-increasing in the homological degree. We wish to understand the maps in this sequence.  To aid this cause, we determine how the ranks of the groups in each homological dimension differ between the cases when $t>T>0$ and $t<-T<0$.
\begin{lemma}
	\label{lem:rank}
	Let $\companion\subset S^3$ be a knot with genus $g$, and let $t>T>0$ be an integer so that Theorem \ref{thm:largen2} holds. Then
	$$\left. \ \begin{array}{ll}
		rk \HFKa_*(\WDnm,1) = rk \HFKa_*(\WDnmneg,1) & {\text if} *\ne0,1 \\
		rk \HFKa_1(\WDnm,1) = rk \HFKa_1(\WDnmneg,1) + t -2\tau(K) &   \\
	rk \HFKa_0(\WDnm,1) = rk \HFKa_0(\WDnmneg,1) -t - 2\tau(K) &  \\
\end{array} \right. $$
	
\end{lemma}

\begin{proof}  The lemma follows from Theorem \ref{thm:largen2}, the definition of $\tau(K)$, and the long exact in homology coming from the short exact sequence of chain complexes,
	$$
\begin{CD}
	0 @>>>	\Filt(\companion,j)
	@>{i}>> \CFa(S^3) @>{p}>> \frac{\CFa(S^3)}{\Filt(\companion,j)} @>>> 0 \end{CD}
.$$

\noindent Since $\HFa(S^3)\cong \F_{(0)}$, the long exact sequence shows that $H_{*-1}(\Filt(\companion,j))\cong H_{*}(\frac{\CFa(S^3)}{\Filt(\companion,j)})$ if $*\ne0,-1$, from which the first part of the lemma follows (taking into account the grading shift in the first part of Theorem \ref{thm:largen2}).  For the second two parts, recall that $\tau(K)$ is defined as:
$$\tau(K)=\mathrm{min}\{j\in\Z|i_*:H_*(\Filt(K,j))\longrightarrow \HFa(S^3)\ {\text{~is non-trivial}}\}.$$

\noindent In the long exact sequence we have:
	$$
\begin{CD}
	0 @>>> H_{1}(\frac{\CFa(S^3)}{\Filt(\companion,j)}) @>{\delta}>>	H_0(\Filt(\companion,j))
	@>{i_*}>> \Z_{(0)} @>{p_*}>> 
	\\ H_{0}(\frac{\CFa(S^3)}{\Filt(\companion,j)}) @>{\delta}>>	H_{-1}(\Filt(\companion,j)) @>>> 0 \end{CD}
,$$
\noindent and the map $i_*$ is trivial precisely when $j < \tau(K)$ and non-trivial otherwise.  Taking the sum over each $j$ from $-g,\ldots ,g$, and examining ranks yields the second two parts of the lemma.

\end{proof}

Next, we observe that the map $f_3$ in the Skein exact sequence, which {\em a priori} is non-increasing in the absolute degree, in fact preserves degree. 

\begin{prop}
	In the exact sequence relating \HFD,\newline \HFDm, $\HFKa(H,1)$,  the map $f_3$ preserves degree. 
 \end{prop}
 \begin{proof}  This follows from the fact that the maps $f_1$ and $f_2$ lower degree by one-half and $f_3$ is non-increasing in the degree, together with the preceding lemma and the fact that the Floer homology of $H$ is supported in degree one-half.  
 \end{proof}

\noindent{\bf Claim:} In the $2t$ applications of the skein sequence connecting \WD \ and \WDneg, $f_2$ is trivial exactly $t-2\tau(K)$ times.

\begin{proof} The preceding proposition tells us that $f_3$ preserves degree. Since $\HFKa(H,1)$ has rank one, supported in degree one-half, \HFDm \ is determined by \HFD \ and whether or not $f_2$ is trivial.  This holds for all $t$, and Lemma \ref{lem:rank} tells us the difference in the ranks of the groups for $t>T>0$ and $-t<-T<0$, thus determining the number of times $f_2$ is trivial.
\end{proof}

Next, recall Proposition $2.4$ of \cite{STau}:

\begin{prop} 
	\label{prop:skein}
 In the exact sequence above, the map $f_2$ is non-trivial if and only if  $\tau(\WDnmm)=1$.  Furthermore, if $\tau(\WDnmm)\ne1$, then it is equal to $0$.
\end{prop}

It is proved in \cite{Livingston1} and \cite{FourBall} that $\tau(K)$ satisfies the following inequality under the operation of changing a crossing in a projection of $K$: $$\tau(K_+)-1\le \tau(K_-) \le \tau(K_+),$$

\noindent where $K_+$ (resp. $K_-$) denote the diagram with the positive (resp. negative) crossing. 
Since each application of the skein sequence arose from changing a single negative crossing to a positive crossing, the above inequality becomes (for $k>0$):

$$\tau(D_+(\companion,t'-k))-k\le \tau(D_+(\companion,t')) \le \tau(D_+(\companion,t'-k)).$$

If $f_2$ were non-trivial for some $t'$ and trivial for $t'-k$, then Proposition \ref{prop:skein} would imply $\tau(D_+(\companion,t'-1))=1$ and $\tau(D_+(\companion,t'-k-1))= 0$, violating the inequality.  Thus $f_2$ is trivial for the first $t-2\tau(K)$ applications of the skein sequence and non-trivial thereafter.   Now \HFDm \ is determined by \HFD \ and knowledge of $f_2$, while $\tau(\WDnm)$ is determined by $f_2$, so this completes our proof of the proposition.\end{proof}

\section{Computation of $\widehat{HFK}(D_+(K,t),0)$ and higher differentials}
\label{sec:middle}
In this section, we complete the calculation of the filtered chain homotopy type of $\HFKa(\WDnm)$.  To this end, recall that the knot Floer homology groups can themselves be thought of as a filtered chain complex, endowed with a differential which strictly lowers the filtration grading.  The homology of $\HFKa(\WDnm)$ under this differential is $\HFa(S^3)\cong\Z_{(0)}$.  The differential is composed of three distinct homomorphisms:
$$d_2 : \HFKa_*(\WDnm,1) \longrightarrow \HFKa_{*-1}(\WDnm,-1)$$
$$d_1^i: \HFKa_*(\WDnm,i) \longrightarrow \HFKa_{*-1}(\WDnm,i-1), \ i=1,0.$$
\noindent Furthermore, the maps $d_1^i$ are induced by chain maps $$\partial_1^i: \CFKa_*(\WDnm,i)\longrightarrow \CFKa_{*-1}(\WDnm,i-1),$$
\noindent defined by counting holomorphic disks which satisfy $n_z(\phi)=1,n_w(\phi)=0$.  

The following is a useful algebraic lemma for the case at hand:
\begin{lemma}
	\label{lemma:d2} Suppose $g(K)=1$ and $\tau(K)=0$ (resp $\tau(K)=1$).  Then the following are equivalent (up to filtered chain homotopy equivalence):
	\begin{enumerate}
		\item $d_2=0$
		\item $d_1^0$ is surjective.
		\item $\HFKa_*(K,0)\cong \Z_{(0)}\oplus \HFKa_{*+1}(K,1)\oplus \HFKa_{*-1}(K,-1)$ \newline(resp. $\cong  \HFKa_{*+1}(K,1)\oplus \HFKa_{*-1}(K,-1)/\Z_{(0)}$)
	\end{enumerate}
\end{lemma}

\begin{proof}
	Each equivalence is a straightforward consequence of the fact that the homology of the chain complex $(\HFKa(K),d_2+d_1^1+d_1^0)$ is isomorphic to $\Z$, supported in grading $0$, with $\tau(K)$ equal to the filtration grading of the generator of this homology.  To see, for example, that $(2)$ implies $(1)$, assume that $d_1^0$ is surjective and rk(Im $d_2)>0$.  Then we can form cycles by adding the chains in $\HFKa(K,1)$ which map non-trivially under $d_2$ to chains in $\HFKa(K,0)$ which map to the same elements in $\HFKa(K,-1)$.  If rk(Im $d_2)>1$, this contradicts the fact that rk($H_*(\HFKa(K),d_2+d_1^1+d_1^0))=1$, while if rk(Im $d_2)=1$ it is easy to see that by a filtered change of basis, the resulting chain complex is filtered chain homotopy equivalent to one with rk(Im $d_2)=0$. The rest of the implications follow by similar considerations and are left to the reader.\end{proof}

\begin{figure}
\begin{center}
	\psfrag{Front}{Front}
\psfrag{m}{$\beta_2=\mu$}
\psfrag{l}{$\lambda_P$}
\psfrag{mk}{$\mu_P$}
\psfrag{b1}{$\beta_1$}
\psfrag{w}{$z$}
\psfrag{z}{$w$}
\psfrag{y}{$y$}
\psfrag{x}{$x$}
\psfrag{a1}{$a_1$}
\psfrag{a2}{$a_2$}
\psfrag{a3}{$a_3$}
\psfrag{a4}{$a_4$}
\psfrag{A}{$\Dom_A$}
\psfrag{B}{$\Dom_B$}
\psfrag{0}{0}
\psfrag{alpha1}{$\alpha_1$}
\caption{\label{fig:Subcomplex}
Depiction of possible domains of Whitney disks, $\phi$, connecting a generator $\{a_2\}\times \p$ to generators $\{a_4\}\times \q$ or $\{y\}\times \q^\prime$, restricted to the Heegaard surface for the pattern.  Since $\phi$ is oriented to go {\em from} $a_2$, the boundary of the domain of $\phi$ must be oriented as shown by the arrows in the figure.  The requirement that $\phi$ satisfies $n_z(\phi)=n_w(\phi)=0$ implies the domain of $\phi$ must have non-zero multiplicity in one, or both, of the domains, $\Dom_A, \Dom_B$.  The orientation of the boundary of $\phi$ and the inward normal orientation of the Heegaard surface imply the multiplicity of $\phi$ in $\Dom_A$ or $\Dom_B$ is negative.}
\includegraphics{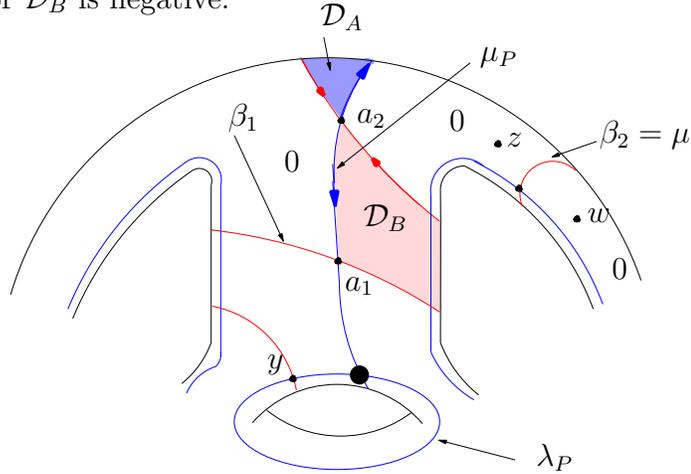}
\end{center}
\end{figure}

	We will show that $(2)$ holds, and hence that $(3)$ (and Theorem \ref{thm:main}) hold.  Recall from Section \ref{sec:meridian}, that the chain complex for the middle filtration level took the following form:
$$\CFKa(\WDnm,0)\sim [\{a_2\}\times \CFa(\Surgnm)] \oplus [\{a_4\}\times \CFa(\Surgnm)] \oplus [\{y\}\times \CFKa(C)].$$
\noindent  We first prove the following:
\begin{prop}
Under the splitting of the generators of $\CFKa(\WDnm,0)$ given above, the generators $[\{a_2\}\times \CFa(\Surgnm)]$ form a subcomplex.  
\end{prop}
\begin{proof}
	To see this, we only have to show that the boundary of any generator of the form $\{a_2\}\times \p$ consists of generators of the same form.  Assume otherwise, that there exists a holomorphic Whitney disk, $\phi$ connecting $\{a_2\}\times \p$ to a generator of the form $\{a_4\}\times \q$ or $\{y\}\times \q^\prime$ satisfying $n_z(\phi)=n_w(\phi)=0$.  In order for $\phi$ to satisfy $n_z(\phi)=n_w(\phi)=0$ and simultaneously be a disk oriented {\em from} $\{a_2\}\times\q$, it must have negative multiplicity in one or both of the regions illustrated in Figure \ref{fig:Subcomplex}. This contradicts the fact that $\phi$ is holomorphic.
\end{proof}

Since $[\{a_2\}\times \CFa(\Surgnm)]$ is a subcomplex, the restriction of the chain map $\partial_1^0$ will be chain map:
$$(\partial_1^0)_{[\{a_2\}\times \CFa(\Surgnm)]}: [\{a_2\}\times \CFa(\Surgnm)] \longrightarrow \CFKa(\WDnm,-1).$$
\noindent Our next claim is the following:
\begin{prop}
	The restriction of $\partial_1^0$ to $[\{a_2\}\times \CFa(\Surgnm)]$ induces an isomorphism on homology.
\end{prop}
\begin{proof}
	There is a canonical ``small'' Whitney disk connecting $\{a_2\}\times \p$ to $\{a_3\}\times \p$ for any $(g-1)$-tuple, $\p$, and which satisfies $n_z(\phi)=1,n_w(\phi)=0$. See Figure \ref{fig:Disk1}. The domain of this disk is topologically an annulus, and can be seen to admit a unique holomorphic representative for a suitably generic choice of almost complex structure on Sym$^g(\Sigma_g)$, see \cite{HolDisk,HolDiskTwo}.
	In the standard way (see \cite{Knots,HolDiskTwo}), we can filter the chain map $\partial_1^0$ with respect to negative area of domains of disks.  With respect to this filtration, the restricted chain map is an isomorphism induced by the aforementioned small disks plus lower order terms, and hence induces an isomorphism on homology. 
\end{proof}

This completes our proof of Theorem \ref{thm:main}: the restricted chain map can be factored as $\partial_1^0\circ i$, where $i$ is the inclusion:
$$i: [\{a_2\}\times \CFa(\Surgnm)] \longrightarrow \CFKa(\WDnm,0).$$ \noindent   Since the restriction induces an isomorphism on homology by the above proposition, we see that the map induced by $\partial_1^0$ (which is $d_1^0$) is surjective.

\section{Preliminary Applications}
\label{sec:examples}

We conclude with two simple applications of Theorem \ref{thm:main} and some qualitative remarks. 

\subsection{Iterated Doubles of the figure eight knot}
We have the following closed formula for the Floer homology of the iterated untwisted Whitehead doubles of the figure eight knot:

\begin{prop}
	\label{prop:fig8}
	Let $4_1$ be the figure eight knot and let $D^n$ denote the $n$-th iterated untwisted double of $4_1$ i.e. $D^0=4_1, D^n=D_+(D^{n-1},0)$ Then we have:
	
	$$\HFKa_*(D^n,i)\cong \left\{\begin{array}{ll} 
		\ \ \ \ \  \bigoplus_{k=0}^{n} \Z_{(1-k)}^{2^n\binom{n}{k}} & i=1\\
		 \Z_{(0)} \bigoplus_{k=0}^{n} \Z_{(-k)}^{2^{n+1}\binom{n}{k}} & i=0\\
		\ \ \ \ \  \bigoplus_{k=0}^{n} \Z_{(-1-k)}^{2^n\binom{n}{k}} & i=-1\\

		0 & otherwise\\
	\end{array}\right.$$
Furthermore, the induced differential $d_1^1$ is injective, $d_1^0$ is surjective, and $d_2$ is zero. Hence $\tau(D^n)=0$
\end{prop}
\begin{proof}
		In order to apply Theorem \ref{thm:main} iteratively, we use our knowledge of the induced differentials $d_1^1,d_1^0,d_2$ acting on $\HFKa(\WDnm)$ which was determined in the preceding section.  We proceed by induction on $n$.  For $n=0$, we have $D^0=4_1=D_+(U,1)$.  $\HFKa(4_1)$ can be determined using various methods (see, for example \cite{Ras1,AltKnots}) but we choose to use Theorem \ref{thm:main}.  Recall that the unknot has Floer homology
		$$\HFKa_*(U,0) \cong \Z_{(0)},$$
		and that $\HFKa_*(U,i)\cong 0$ for $i\ne0$.  
		This immediately implies that $H_*(\Filt(U,0))\cong \Z_{(0)}$ and that $\tau(U)=0$.  Hence we see that $\bigoplus_{i=-0}^{0}[H_{*}(\Filt(U,i))] \cong Z_{(0)}$.  Plugging this result into the formula of Theorem \ref{thm:main} for $\HFKa(D_+(U,1))$ (where we use the parameters $t=1$, $\tau(U)=0$, and $g(U)=0$), we see that 
	$$\HFKa_*(D_+(U,1),i)\cong \left\{\begin{array}{ll} 
		\Z_{(1)}^{-1} \oplus \Z_{(1)}^2 = \Z_{(1)} & i=1\\
		\Z_{(0)}^{-1} \oplus \Z_{(0)}^4 = \Z_{(0)}^3 & i=0\\
		\Z_{(-1)}^{-1} \oplus \Z_{(-1)}^2 = \Z_{(-1)} & i=-1\\

		0 & otherwise\\
	\end{array}\right.$$
Furthermore, we know from the previous section and the fact that we are in the case when $t=1>2\tau(U)$ that $d_1^1$ is injective, $d_1^0$ is surjective, and $d_2$ is zero.  This completes the base case.

Assume that the proposition holds for $n$.  This implies that
$$H_*(\Filt(D^n,-1))) \cong \HFKa(D^n,-1) \cong  \bigoplus_{k=0}^{n} \Z_{(-1-k)}^{2^n\binom{n}{k}}.$$
As for $H_*(\Filt(D^n,0)))$ we have,
$$H_*(\Filt(D^n,0))) \cong H_*( \HFKa(D^n,-1)\oplus \HFKa(D^n,0), d_1^0)\cong $$ 
$$H_*(\bigoplus_{k=0}^{n} \Z_{(-1-k)}^{2^n\binom{n}{k}} \oplus	 \Z_{(0)} \bigoplus_{k=0}^{n} \Z_{(-k)}^{2^{n+1}\binom{n}{k}}, d_1^0)\cong 	 \Z_{(0)} \bigoplus_{k=0}^{n} \Z_{(-k)}^{2^{n}\binom{n}{k}},$$

Where the first congruence is the definition of $\Filt(K,0)$ for a genus one knot, the second follows from our inductive hypothesis, and the final follows from the fact that $d_1^0$ is assumed to be surjective.

For the final filtration, we clearly have $H_*(\Filt(D^n,1))\cong \Z_{(0)}.$
Thus, we see that
$$\bigoplus_{i=-1}^{1}[H_{*}(\Filt(D^n,i))]  \cong [\bigoplus_{k=0}^{n} \Z_{(-1-k)}^{2^n\binom{n}{k}}] \oplus [ \Z_{(0)} \bigoplus_{k=0}^{n} \Z_{(-k)}^{2^{n}\binom{n}{k}}] \oplus [Z_{(0)}] = \Z_{(0)}^2 \bigoplus_{k=0}^{n+1} \Z_{(-k)}^{2^n\binom{n+1}{k}} $$

Applying Theorem \ref{thm:main} with parameters $t=0$, $\tau(D^n)=0$, $g(D^n)=1$, we have that:
$$\HFKa_*(D^{n+1},i)\cong \left\{\begin{array}{ll} 
	\Z_{(1)}^{-4}	\bigoplus_{i=-1}^{1}[H_{*-1}(\Filt(D^n,i))]^2 = 	\Z_{(1)}^{-4} \oplus \Z_{(1)}^4 \bigoplus_{k=0}^{n+1} \Z_{(-k)}^{2^{n+1}\binom{n+1}{k}} & i=1\\
	\Z_{(0)}^{-7}	\bigoplus_{i=-1}^{1}[H_{*}(\Filt(D^n,i))]^4 =	\Z_{(0)}^{-7} \oplus \Z_{(0)}^8 \bigoplus_{k=0}^{n+1} \Z_{(-k)}^{2^{n+2}\binom{n+1}{k}}& i=0\\
	\Z_{(-1)}^{-4} \bigoplus_{i=-1}^{1}[H_{*+1}(\Filt(D^n,i))]^2 = 	\Z_{(-1)}^{-4} \oplus \Z_{(-1)}^4 \bigoplus_{k=0}^{n+1} \Z_{(-k)}^{2^{n+1}\binom{n+1}{k}} & i=-1.\\

		0 & otherwise\\
	\end{array}\right.$$
	\noindent Quotienting by the negative exponents in the above equation yields the formula given by Proposition \ref{prop:fig8} for $D^{n+1}$, thus completing the inductive step.  Theorem \ref{thm:main} and Lemma \ref{lemma:d2} show that $d_1^1$ is injective, $d_1^0$ is surjective, and $d_2=0$.

\end{proof}
 We found it notable that while the figure eight knot is an alternating knot and has particularly simple Floer homology, by forming its iterated untwisted doubles we obtain knot Floer homology groups which become incredibly complicated.  In particular, the {\em width} of the Floer homology (the number of diagonals on which knot Floer homology - plotted on a grid whose axes are the homological and filtration grading - is supported) can be made arbitrarily large.  Indeed, the width grows linearly with the number of times we double.  Also, the total rank of the Floer homology grows exponentially with the number of times we double.  In some sense, the Floer homology of the Whitehead double is ``seeing'' all the Floer homology of the companion knot.

	\subsection{Surgery on Whitehead doubles}
	As a final application, we determine the Floer homology of $+1$-surgery on the Whitehead double of a knot, $K$.
	\begin{prop}
		Let $S^3_{+1}(\WDnm)$ denote the manifold obtained by $+1$-surgery on \WD.  Then for $t\ge 2\tau(K)$ we have:
		$$\HFa_*(S^3_{+1}(\WDnm)) \cong \Z_{(-1)}^{t-2g-2}\oplus \Z_{(0)}^{t-2g-1}\bigoplus_{i=-g}^{g}[H_{*+1}(\Filt(K,i))]^2  \bigoplus_{i=-g}^{g}[H_{*}(\Filt(K,i))]^2 $$
		While for $t<2\tau(K)$ we have:
			$$\HFa_*(S^3_{+1}(\WDnm)) \cong $$ $$\Z_{(-1)}^{4\tau(K)-t-2g-2}\oplus\Z_{(-2)}^{2\tau(K)-t}\oplus \Z_{(0)}^{2\tau(K)-2g-2}\bigoplus_{i=-g}^{g}[H_{*+1}(\Filt(K,i))]^2  \bigoplus_{i=-g}^{g}[H_{*}(\Filt(K,i))]^2 $$
	\end{prop}
	\begin{proof} This is a straightforward application of Theorem $4.4$ of \cite{Knots}, together with Theorem \ref{thm:main}.  For a genus one knot, $K$, Theorem $4.4$ of \cite{Knots} identifies 
		$$\HFa_*(S^3_{+1}(K)) \cong H_{*}(C\{\mathrm{max}(i,j)=0\}). $$
\noindent And this latter group is equal to:
$$H_*( \HFKa_*(K,1)\{-2\}\oplus \HFKa_*(K,0) \oplus \HFKa_*(K,-1), d_1^0 + \tilde{d}_1^0),$$		
\noindent where $d_1^0: \HFKa_*(K,0)\rightarrow \HFKa_{*-1}(K,-1)$ is the map induced by the chain map $\partial_1^0$ discussed in Section \ref{sec:middle}. The map $\tilde{d}_1^0: \HFKa_*(K,0)\rightarrow \HFKa_{*+1}(K,1)$ is induced by the chain map $\tilde{\partial}_1^0: \HFKa_*(K,0)\rightarrow \HFKa_{*+1}(K,1)$ which counts pseudo-holomorphic Whitney disks satisfying $n_w(\phi)=1,n_z(\phi)=0$.  The $\{-2\}$ indicates that we shift the grading of $\HFKa_*(K,1)$ down by $2$ (this is induced by the action of $U$ on $\CFKinf(K)$).  In the present situation, Theorem \ref{thm:main} informs us of the groups $\HFKa_*(\WDnm,1)[-2]\oplus \HFKa_*(\WDnm,0) \oplus \HFKa_*(\WDnm,-1)$, and furthermore that $d_1^0$ is surjective.  It follows from algebraic properties of $\CFKinf(K)$ that  $\tilde{d}_1^0$ will also be surjective.  Alternatively, this can be seen by the same method used in Section \ref{sec:middle} to show that $d_1^0$ is surjective.  The proposition follows.
\end{proof}

\commentable{
\bibliographystyle{plain}
\bibliography{biblio}
}

\end{document}